\documentclass[11pt]{amsart}

\usepackage{amsmath,amssymb,amsthm,dsfont}
\usepackage[colorlinks,linkcolor=blue,anchorcolor=blue,citecolor=blue]{hyperref}

\newcommand{\NN}{\mathbb{N}}
\newcommand{\ZZ}{\mathbb{Z}}
\newcommand{\RR}{\mathbb{R}}
\newcommand{\CC}{\mathbb{C}}
\newcommand{\cX}{\mathcal{X}}
\newcommand{\cR}{\mathcal{R}}
\newcommand{\cA}{\mathcal{A}}
\DeclareMathOperator{\1}{\mathds{1}}
\newcommand{\Mod}[1]{\ (\operatorname{mod}#1)}

\theoremstyle{plain}
\newtheorem{theorem}{Theorem}[section]
\newtheorem{corollary}[theorem]{Corollary}
\newtheorem{lemma}[theorem]{Lemma}
\newtheorem{proposition}[theorem]{Proposition}
\newtheorem{remark}{Remark}[section]

\numberwithin{equation}{section}

\allowdisplaybreaks

\title[Poissonian pair correlation and energy estimates]{Metric Poissonian pair correlation for real sequences and energy estimates}

\author[B. Kerr]{Bryce Kerr}
\address{B.K.: School of Science, University of New South Wales, Canberra, ACT, Australia}
\email{bryce.kerr89@gmail.com}

\author[H. Wang]{Hongliang Wang}
\address{H.W.: School of Science, University of New South Wales, Canberra, ACT, Australia}
\email{hongliang.wang00@gmail.com}

\begin{document}

\begin{abstract}
We establish new conditions under which a sequence of real numbers has metric Poissonian pair correlation. These conditions strengthen results of Aistleitner, El-Baz and Munsch (2021) and resolve one of their open problems under a mild growth assumption. As applications, we show that quantitatively convex and polynomial sequences have metric Poissonian pair correlation. 
\end{abstract}

\maketitle

\tableofcontents

\section{Introduction}

A sequence of distinct real numbers $(x_n)_{n\in\NN}$, which we always assume to be in ascending order, is said to have Poissonian pair correlation if for any $s>0$, 
    \begin{equation}
    \label{eq:Poissonian pair correlation}
    \lim_{N\to\infty}\frac{1}{N}\left|\left\{1\leq n_1\neq n_2\leq N:\|x_{n_1}-x_{n_2}\|\le \frac{s}{N}\right\}\right|= 2s, 
    \end{equation}
where $\|\cdot\|$ denotes the distance to the nearest integer. 

Poissonian pair correlation is a fine-scale statistic for the distribution of fractional parts and can be considered as capturing pseudo-randomness, since a sequence of independent, uniformly distributed random variables on $[0,1)$ satisfies \eqref{eq:Poissonian pair correlation} with probability $1$. For a statistical perspective on this notion and its connection to uniform distribution, see~\cite{ALL18,LS20}. 

The study of pair correlation of fractional parts was initiated by Rudnick and Sarnak~\cite{RS98}, who proved that for any integer $d\geq 2$, there exists a set $P\subseteq\RR$ of full Lebesgue measure such that for any $\alpha\in P$, the sequence $(\alpha n^{d})_{n\in\NN}$ has Poissonian pair correlation. This problem has origins in quantum physics, where the cases $d=1,2$ arise naturally in the analysis of energy levels of physical systems, as investigated by Berry and Tabor~\cite{BT77} (see also~\cite{M01}). The famous Berry--Tabor conjecture predicts that the local spacing of energy levels of classically integrable systems admits the same limiting distribution as that of a Poisson point process. Such a distribution is known to be determined by correlations of all levels, with pair correlation  being the most tractable to analyse, see for example~\cite{RSZ01}. 

Very few specific sequences are known to have Poissonian pair correlation. A notable recent result by Lutsko, Sourmelidis and Technau~\cite{LST24} shows that the sequence $(\alpha n^\theta)_{n\in\NN}$ has Poissonian pair correlation for all $\alpha>0$ and $0<\theta<14/41$, and this was subsequently extended to $0<\theta<43/117$ by Radziwi{\l}{\l} and Shubin~\cite{RS24}. We refer the reader to these papers for further discussion of this type of problems. 

Results concerning sequences of the form $(\alpha x_n)_{n\in\NN}$ for almost all $\alpha\in\RR$, in the sense of full Lebesgue measure, belong to the metric theory. When $(\alpha x_n)_{n\in\NN}$ has Poissonian pair correlation for almost all $\alpha$, we say that the sequence $(x_n)_{n\in\NN}$ is metric Poissonian. As previously mentioned, the sequence $(n^d)_{n\in\NN}$ is metric Poissonian. This result was generalised by Boca and Zaharescru~\cite{BZ00} to polynomial sequences $(p(n))_{n\in\NN}$, where $p\in\ZZ[X]$ has degree at least $2$. See~\cite{BE25+,CLZ15,RZ99,W18} for further metric results involving integer-valued and rational-valued sequences. 

Over the past decade, significant progress has been made in establishing conditions that determine when an integer valued sequence is metric Poissonian. Let $E_N$ denote the truncated additive energy of the sequence $(x_n)_{n\in\NN}$, which is defined by 
    \[
    E_N = |\{1\le n_1,n_2,n_3,n_4\leq N:x_{n_1}+x_{n_2}=x_{n_3}+x_{n_4}\}|. 
    \]
A foundational result of Aistleitner, Larcher and Lewko~\cite[Theorem~1]{ALL17} is as follows. 

\begin{theorem}
\label{thm:ALL}
Let $(x_n)_{n\in\NN}$ be a sequence of distinct integers. Assume that there exists $\delta>0$ such that 
    \begin{equation}
    \label{eq:additive energy trivial bound}
    E_N\ll N^{3-\delta}.
    \end{equation}
Then $(x_n)_{n\in\NN}$ is metric Poissonian. 
\end{theorem}

The exponent $3-\delta$ here is sharp, since trivially we have 
    \[
    E_N\le N^3, 
    \]
and Bourgain~\cite[Appendix]{ALL18} has shown that $(x_n)_{n\in\NN}$ is not metric Poissonian provided
    \begin{equation}
    \label{eq:not metric Poissonian}
    \limsup_{N\rightarrow \infty}\frac{E_N}{N^3}>0. 
    \end{equation}
There is a gap between conditions~\eqref{eq:additive energy trivial bound} and~\eqref{eq:not metric Poissonian}, and it remains an open problem whether there exists an additive or combinatorial characterisation of the metric Poissonian property. It follows from work of Aistleitner, Lachmann and Technau~\cite{ALT19} that a criterion based solely on additive energy is insufficient. We refer the reader to~\cite{BE25+,BCGW18,BW20,LS20b} for progress on related extremal problems. 

The question of whether some variation of Theorem~\ref{thm:ALL} holds for real-valued sequences is more difficult. Progress on this problem was made by Rudnick and Technau~\cite{RT20}, who reduced to a lattice point counting problem, and subsequently Aistleitner, El-Baz and Munsch~\cite{AEM21} showed that this lattice point counting problem can be controlled by an approximate additive energy $E^*_N$, defined by 
    \[
    E^*_N = |\{1\leq n_1,n_2,n_3,n_4\leq N:|x_{n_1}-x_{n_2}+x_{n_3}-x_{n_4}|<1\}|. 
    \]
The current state of the art is given by~\cite[Theorem~1]{AEM21}. 

\begin{theorem}
\label{thm:AEM}
Let $(x_n)_{n\in\NN}$ be a sequence of positive real numbers for which there exists $c>0$ such that 
    \[
    x_{n+1}-x_n \geq c. 
    \]
Assume that there exists $\delta>0$ such that 
    \[
    E^{*}_N \ll N^{183/76-\delta}. 
    \]
Then $(x_n)_{n\in\NN}$ is metric Poissonian. 
\end{theorem}

As an application of their techniques, the authors of~\cite{AEM21} showed that the sequence $(n^\theta)_{n\in\NN}$ is metric Poissonian for $\theta>1$. A related result was established by Rudnick and Technau~\cite{RT22}, who proved that $(n^\theta)_{n\in\NN}$ is metric Poissonian for $0<\theta<1$. See also~\cite{EMV15,RT20} for further metric results involving real-valued sequences. 

The main objective of this paper is to further develop the understanding of metric Poissonian pair correlation for real-valued sequences. 

\section{Main results}

For $0<\gamma\leq1$, define $E^*_{N,\gamma}$ to be the number of solutions to the inequality 
    \[
    |x_{n_1}-x_{n_2}+x_{n_3}-x_{n_4}| < \gamma 
    \]
for variables $1\leq n_1,n_2,n_3,n_4\leq N$. When $\gamma=1$, we write 
    \[
    E^*_{N} = E^*_{N,1}. 
    \]
Our main result is as follows. 

\begin{theorem}
\label{thm:main}
Let $(x_n)_{n\in\NN}$ be a sequence of positive real numbers for which there exist $c>0$ and $\eta\in(2/3,1]$ such that 
    \begin{equation}
    \label{eq:main spacing condition}
    x_{n+1}-x_n \geq \frac{c}{n^{1-\eta}}. 
    \end{equation}
Assume that there exist $\delta>0$ and $\kappa\in[0,1)$ such that 
    \begin{equation}
    \label{eq:main energy condition}
    E^*_{N,1/N}\ll N^{3-\kappa-\delta} 
    \end{equation}
and
    \begin{equation}
    \label{eq:main growth condition}
    \sum_{\substack{1\leq n_1<n_2\leq N\\x_{n_2}-x_{n_1}\geq1}}\frac{1}{\sqrt{x_{n_2}-x_{n_1}}} \ll N^{1+\kappa/2}. 
    \end{equation}
Then $(x_n)_{n\in\NN}$ is metric Poissonian. 
\end{theorem}

The spacing condition~\eqref{eq:main spacing condition} is significantly weaker than that in Theorem~\ref{thm:AEM} of Aistleitner, El-Baz and Munsch, as it permits sequences whose gaps tend to $0$. We note that the threshold $2/3$ in~\eqref{eq:main spacing condition} arises from a technical aspect of our approach and is likely not sharp.

We expect the condition~\eqref{eq:main growth condition} may be weakened and we refer the reader to the discussion following Proposition~\ref{prop:decreasing property} for an indication of how it appears in our proof.

We next deduce a number of consequences from Theorem~\ref{thm:main}. Our first Corollary is a direct improvement on Theorem~\ref{thm:AEM}, which follows from estimating the left-hand side of~\eqref{eq:main growth condition} in terms of $E^{*}_N$. 

\begin{corollary}
\label{cor:energy only}
Let $(x_n)_{n\in\NN}$ be a sequence of positive real numbers for which there exists $c>0$ such that $x_{n+1}-x_n\geq c$. Assume that there exists $\delta>0$ such that 
    \begin{equation}
    \label{eq:energy only condition}
    E^*_N \ll N^{5/2-\delta}. 
    \end{equation}
Then $( x_n)_{n\in\NN}$ is metric Poissonian. 
\end{corollary}

Comparing the energy bound~\eqref{eq:energy only condition} with that in Theorem~\ref{thm:AEM}, our improvement becomes clear upon noting that 
    \[
    \frac{5}{2} = \frac{183}{76}+\frac{7}{76}.
    \]

We next apply Theorem~\ref{thm:main} to address some open problems posed by Aistleitner, El-Baz and Munsch~\cite{AEM21}. 

In~\cite[Open Problem~4]{AEM21}, the authors ask whether a condition of the form 
    \[
    E^*_{N,\gamma}\ll\gamma N^{4-\delta} \quad \text{where} \quad \gamma\approx\frac{1}{N} 
    \]
is sufficient to establish the metric Poissonian property, as this would give a unified picture for both integer and real valued sequences. Theorem~\ref{thm:main} (with $\kappa=0$) answers this question in the affirmative, subject to the additional assumption \eqref{eq:main growth condition}; that is, 
    \[
    \sum_{\substack{1\leq n_1<n_2\leq N\\x_{n_2}-x_{n_1}\geq1}}\frac{1}{\sqrt{x_{n_2}-x_{n_1}}} \ll N^{1+o(1)}. 
    \]
Our next corollary states that this assumption is satisfied when the sequence does not concentrate in short intervals. 

\begin{corollary}
\label{cor:short interval}
Let $(x_n)_{n\in\NN}$ be a sequence of positive real numbers satisfying~\eqref{eq:main spacing condition}. Assume that there exists $\delta>0$ such that 
    \[
    E^*_{N,1/N}\ll N^{3-\delta}, 
    \]
and for any $X,H\gg 1$, 
    \begin{equation}
    \label{eq:growth condition short interval}
    |\{n\in\NN:x_n\in[X,X+H]\}| \ll H^{1/2}. 
    \end{equation}
Then $(x_n)_{n\in\NN}$ is metric Poissonian. 
\end{corollary}

The condition~\eqref{eq:growth condition short interval} is satisfied, for example, when $x_n$ are consecutive values of a function whose derivative is positive and grows faster than linearly. 

In~\cite[pg. 486]{AEM21}, the authors discuss the problem of whether convex or polynomial sequences are metric Poissonian. By combining Theorem~\ref{thm:main} with estimates for the additive energy of convex sequences due to Shkredov~\cite{S13} and subsequently improved by Bloom~\cite{B25+}, we answer both of these question in the affirmative. 

\begin{corollary}
\label{cor:convex sequences}
Let $(x_n)_{n\in\NN}$ be a sequence of real numbers satisfying 
    \begin{equation}
    \label{eq:convex condition}
    x_{n+1}-x_n \geq x_{n}-x_{n-1}+cn^{-1/10^4}. 
    \end{equation}
Then $(x_n)_{n\in\NN}$ is metric Poissonian. 
\end{corollary}

The exponent $10^{-4}$ in~\eqref{eq:convex condition} is not sharp and we have made no attempt to optimise it. Since non-linear polynomial sequences are convex, they are consequently metric Poissonian. The proof of Corollary~\ref{cor:convex sequences} relies on advanced techniques from additive combinatorics to establish the relevant energy estimate. We provide an alternative approach by establishing a variation of Hua's inequality to show that polynomial sequences are metric Poissonian. 

\begin{corollary}
\label{cor:polynomial sequences}
Let $p\in\RR[X]$ be a polynomial of degree at least $2$. Then $(p(n))_{n\in\NN}$ is metric Poissonian. 
\end{corollary}

The relevant energy estimates for quantitatively convex and polynomial sequences are provided in Propositions~\ref{prop:convex} and~\ref{prop:polynomial}, respectively. These estimates may be of independent interests beyond the metric Poissonian problem. 

There are other interesting arithmetic sequences to which one may attempt to apply Theorem~\ref{thm:main}. One such example is the sequence $(\sigma_{\beta}(n))_{n\in\NN}$ where 
    \[
    \sigma_{\beta}(n) := \sum_{d|n}d^{\beta}.
    \]
Based on numerical computations for a handful of algebraic irrational $\alpha$'s and rational $\beta$'s, the local spacing distributions of $(\alpha\sigma_{\beta}(n))_{n\in\NN}$ appear to be Poissonian when $\beta\neq1$, and so does the pair correlation. However, we have been unable to establish suitable energy estimates in order to apply Theorem~\ref{thm:main}. For $\beta=1$, truncations of $(\alpha\sigma_1(n))_{n\in\NN}$ seem to contain too few gaps to satisfy Poissonian statistics.

\section*{Acknowledgement}

The authors would like to thank Ze\'{e}v Rudnick for an interesting discussion, for bringing to their attention the problem of determining correlations of sequences $\sigma_{\beta}(n)$ and for showing them some numerical insights. 
The authors would like to thank Christoph Aistleitner, Igor Shparlinski and Andrei Shubin for some comments on a previous version of this manuscript. 

During the preparation of this work, B. Kerr was supported by
ARC Grants DE220100859 and DP230100534. 

\section{Overview and structure of the paper}
\label{sec:over}
Our proof of Theorem~\ref{thm:main} begins with the following result of Rudnick and Technau~\cite[Theorem 1.2]{RT20}. 

\begin{theorem}
\label{thm:RT}
Let $(x_n)_{n\in\NN}$ be a sequence of distinct positive real numbers. Assume that there exits $\delta>0$ such that for sufficiently small $\varepsilon>0$, 
    \[
    \sum_{\substack{1\leq m\leq N^{1+\varepsilon}\\1\leq n_1\neq n_2\leq N\\m|x_{n_1}-x_{n_2}|<N^\varepsilon}}1 \ll N^{2-\delta} 
    \]
and 
    \begin{equation}
    \label{eq:overview RT condition}
    \sum_{\substack{1\leq |m_1|,|m_2|\leq N^{1+\varepsilon}\\1\leq n_1\neq n_2\leq N\\1\leq n_3\neq n_4\leq N\\|m_1(x_{n_1}-x_{n_2})-m_2(x_{n_3}-x_{n_4})|<N^\varepsilon}}1 \ll N^{4-\delta}. 
    \end{equation}
Then $(x_n)_{n\in\NN}$ is metric Poissonian. 
\end{theorem}

Theorem~\ref{thm:RT} reduces Theorem~\ref{thm:main} to a lattice point counting problem over a family of lattices. To establish \eqref{eq:overview RT condition}, it suffices to bound sums of the form 
    \begin{equation}
    \label{eq:S def}
    S(\cX,\alpha,M,K) := \sum_{\substack{1\leq m_1,m_2\leq M \\ x_1,x_2\in\cX \\|m_1x_1-m_2x_2|\leq K}}\alpha(x_1)\alpha(x_2) 
    \end{equation}
where 
    \[
    \cX := \{x_{n_2}-x_{n_1}:1\leq n_1<n_2\leq N\}
    \]
and 
    \[
    \alpha(x) := |\{1\leq n_1<n_2\leq N:x_{n_2}-x_{n_1}=x\}|. 
    \]
Note that $E^*_{N,\gamma}$ is an approximate $\ell_2$-norm of $\alpha$: 
    \[
    E^*_{N,\gamma} \approx \sum_{\substack{x_1,x_2\in\cX\\|x_1-x_2|\leq\gamma}}\alpha(x_1)\alpha(x_2). 
    \]

Proceeding as in~\cite[Section~5]{AEM21}, we estimate~\eqref{eq:S def} in terms of the second moment of a Dirichlet polynomial of the following shape: 
    \begin{equation}
    \label{eq:overview mean value}
    \frac{1}{T}\int_{-T}^{T}\left|\sum_{x\in\cX}\frac{\alpha(x)}{x^{it}}\sum_{N<n\leq2N}\frac{1}{n^{it}}\right|^2dt, 
    \end{equation}
where $T$ depends on the range of the support of $\alpha$. 

In~\cite{AEM21}, the authors first carry out some combinatorial preparations to control the size of $T$. After an application of H\"{o}lder's inequality, this leads to a reduction to a higher order moment of the Riemann zeta function and an approximate $\ell_2$-norm of $\alpha$, which is essentially $E^*_{N,\gamma}$ for some $\gamma$ depending on $T$. This procedure loses information. For example, controlling the size of $T$ restricts the reduction to $E^{*}_N$, as opposed to the finer quantity $E^{*}_{N,1/N}$. Moreover, applying H\"{o}lder's inequality while still reducing to an $\ell_2$ bound for $\alpha$ forces the use of an $\ell_{\infty}$ bound for the Dirichlet polynomial supported on $\alpha$, and for general sequences $(x_n)_{n\in\NN}$, the best available bound is trivial. 

It is at this stage that our approach diverges from~\cite{AEM21}. To estimate the mean values~\eqref{eq:overview mean value}, we employ some ideas of Heath-Brown from~\cite{H79,H18}. Specifically, we will establish a number of inequalities for the sums $S(\cX,\alpha,M,K)$, defined as in~\eqref{eq:S def}, in terms of $S(\cX,\alpha,M',K')$ with parameters $M'$ and $K'$ of different sizes. These inequalities will then be combined to derive a recursive bound for $S(\cX,\alpha,M,K)$.

In Section~\ref{sec:lattice}, we establish our main result for the lattice point counting problem~\eqref{eq:S def} and state the aforementioned recursive inequalities. A proof of the main result, conditional on these inequalities is given here, while the proofs of the inequalities are deferred to Sections~\ref{sec:decomposition},~\ref{sec:geometry} and~\ref{sec:dirichlet}. 

In Section~\ref{sec:poissonian}, we provide details on how our estimate for~\eqref{eq:S def} can be combined with Theorem~\ref{thm:RT} to establish the main result, Theorem~\ref{thm:main}. We also derive Corollaries~\ref{cor:energy only} and~\ref{cor:short interval} by addressing the growth condition~\eqref{eq:main growth condition}. 

Finally, in Section~\ref{sec:energy}, we compute the relevant energy estimates for quantitatively convex and polynomial sequences, thereby establishing Corollaries~\ref{cor:convex sequences} and~\ref{cor:polynomial sequences}. 

\section{Lattice point counting}
\label{sec:lattice}

We establish estimates for the lattice point counting problem~\eqref{eq:S def} for general sets $\cX$ and weight functions $\alpha$ supported on $\cX$. Our bounds are expressed in terms of the $\ell_1$-norm 
    \begin{equation}
    \label{eq:l1 norm def}
    \|\alpha\|_1 := \sum_{x\in\cX}\alpha(x), 
    \end{equation}
the approximate $\ell_2$-norm 
    \begin{equation}
    \label{eq:l2 norm def}
    \|\alpha\|_{2,N} := \left(\sum_{\substack{x_1,x_2\in\cX\\|x_1-x_2|\leq1/N}}\alpha(x_1)\alpha(x_2)\right)^{1/2} 
    \end{equation}
and a weighted sum over $\cX$ corresponding to the left-hand side of \eqref{eq:main growth condition}. 

\begin{theorem}
\label{thm:dirichlet}
Let $M,K\gg1$, $\cX\subseteq[1,\infty)$ be a finite set, and $\alpha$ be a function on $\cX$ taking positive real values. Then
    \begin{align*}
    S(\cX,\alpha,M,K) &\ll M^{3/2+o(1)}K^{1+o(1)}\|\alpha\|_{2,M}\sum_{x\in\cX}\frac{\alpha(x)}{x^{1/2}} \\
    &\quad +(MK)^{1/2+o(1)}\|\alpha\|_1\|\alpha\|_{2,M}+(MK)^{1+o(1)}\|\alpha\|_{2,M}^2. 
    \end{align*}
\end{theorem}

\subsection{Preliminary setup}

We begin by setting up some notation that will be used throughout. We write $n\sim N$ to denote $N<n\leq2N$. Let $\beta=(\beta(n))_{n\in\NN}$ denote a sequence of complex numbers, and define $|\beta|=(|\beta(n)|)_{n\in\NN}$. For $N\gg1$, define 
    \[
    \widetilde{S}(X,\alpha,\beta,N,K) 
    = \sum_{\substack{n_1,n_2\sim N\\x_1,x_2\in X\\
    |n_1x_1-n_2x_2|\leq K}}
    \alpha(x_1)\alpha(x_2)\beta(n_1)\beta(n_2).
    \]
This generalised form involving $\beta$ will be useful when handling mean values of Dirichlet polynomials in Section~\ref{sec:dirichlet}. In most cases, we take $\beta(n)=1$, which leads to the simplified notation 
    \begin{equation}
    \label{eq:dyadic S def}
    \widetilde{S}(X,\alpha,N,K) := \sum_{\substack{n_1,n_2\sim N\\x_1,x_2\in X\\|n_1x_1-n_2x_2|\leq K}}\alpha(x_1)\alpha(x_2). 
    \end{equation}
This differs from $S(\cX,\alpha,M,K)$ only in that the variables $n_1$ and $n_2$ are restricted to a dyadic interval. 

For a non-negative integer $k$, define 
    \[
    \cX_k = \cX\cap [2^k,2^{k+1}). 
    \]
Let $\alpha_k$ denote the restriction of $\alpha$ to $\cX_k$; that is, 
    \[
    \alpha_k(x) = \begin{cases}
    \alpha(x) &$if $ x\in \cX_k, \\ 0 &$otherwise$. 
    \end{cases}
    \]
We also define $\|\alpha_k\|_1$ and $\|\alpha_k\|_{2,N}$ by restricting the summations to $\cX_k$: 
    \[
    \|\alpha_k\|_1 = \sum_{x\in\cX_k}\alpha(x)
    \]
and 
    \[
    \|\alpha_k\|_{2,N} = \left(\sum_{\substack{x_1,x_2\in\cX_k\\|x_1-x_2|\leq1/N}}\alpha(x_1)\alpha(x_2)\right)^{1/2}.
    \]

It is straightforward to show that 
    \begin{equation}
    \label{eq:norm sums}
    \sum_{k=0}^\infty\|\alpha_k\|_1^2 \leq \|\alpha\|_1^2 
    \quad \text{and} \quad 
    \sum_{k=0}^\infty\|\alpha_k\|_{2,N}^2 
    \leq \|\alpha\|_{2,N}^2. 
    \end{equation}
The approximate $\ell_2$-norm $\|\alpha_k\|_{2,N}$ satisfies the following monotonicity property in parameter $N$. 

\begin{lemma}
\label{lem:l2 norm increasing}
For any $1\ll N\ll M$, 
    \[
    \|\alpha_k\|^2_{2,N} \ll \frac{M}{N}\|\alpha_k\|^2_{2,M}. 
    \]
\end{lemma}

\begin{proof}
Applying Lemma~\ref{lem:mean value theorem} below, we obtain 
    \[
    \|\alpha_k\|^2_{2,N} \ll \frac{1}{N}\int_{-N}^N\left|\sum_{x\in\cX_k}\alpha(x)e(xt)\right|^2dt \leq \frac{M}{N}\cdot\frac{1}{M}\int_{-M}^M\left|\sum_{x\in \cX_k}\alpha(x)e(xt)\right|^2dt. 
    \]
The result follows since another application of Lemma~\ref{lem:mean value theorem} yields 
    \[
    \frac{1}{M}\int_{-M}^M\left|\sum_{x\in \cX_k}\alpha(x)e(xt)\right|^2dt \ll \|\alpha_k\|_{2,M}^2. 
    \qedhere
    \]
\end{proof}

\subsection{Statement of recursive inequalities}

In the following statements, we consider $M,N,K\gg1$, $\cX\subseteq[1,\infty)$ a finite set, $\alpha$ a function on $\cX$ taking positive real values, and $\beta$ a sequence of complex numbers. 

Our first proposition is based on a combinatorial decomposition of the summation~\eqref{eq:S def}, which partitions the variables into dyadic intervals and reduces the problem to summations of the form~\eqref{eq:dyadic S def}. 

\begin{proposition}[Dyadic partition]
\label{prop:dyadic partition}
There exists $N\leq 2M$ such that 
    \[
    S(\cX,\alpha,M,K) \ll (MK)^{o(1)}\sum_{k=0}^{\infty}\widetilde{S}(\cX_k,\alpha,N,4K). 
    \]
\end{proposition}

We establish several recursive bounds for $\widetilde{S}$ that result in an inflation of the parameter $K$. To control this, we appeal to the following proposition, which can be viewed as a form of linearity in parameter $K$. 

\begin{proposition}[Linear property]
\label{prop:linear property}
For any $J\geq1$, 
    \[
    \widetilde{S}(\cX_k,\alpha,\beta,N,JK)\ll J\widetilde{S}(\cX_k,\alpha,|\beta|,N,K).
    \]
\end{proposition}

Our first recursive bound for $\widetilde{S}$ draws on the geometry of numbers and decreases the length of summation. It is inspired by a similar argument of Heath-Brown~\cite[Lemma~6]{H18}. 

\begin{proposition}[Decreasing property]
\label{prop:decreasing property}
For any $k\geq0$ and $N_0\leq N$, there exists $L\leq4N_0$ such that 
    \[
    \widetilde{S}(\cX_k,\alpha,N,K) \ll \frac{NK}{2^{k}}\|\alpha_k\|^2_1+\frac{N}{N_0}\|\alpha_k\|^2_1+\frac{N^{1+o(1)}}{N_0}\widetilde{S}\left(\cX_k,\alpha,L,\frac{N_0K}{N}\right). 
    \]
\end{proposition}

The proof of Proposition~\ref{prop:decreasing property} is based on the geometry of numbers and by a volume-packing argument (see, for example,~\cite[Lemma~3.24]{TV10}), we expect a matching lower bound of 
    \[
    \frac{NK}{2^k}\|\alpha_k\|_1^2. 
    \]
This term contributes to the condition~\eqref{eq:main growth condition} of our main theorem. 

The parameter $L$ in Proposition~\ref{prop:decreasing property} arises from an application of the dyadic pigeonhole principle. Our next result allows us to restrict attention to values of $L$ close to the maximal size, $L\approx N_0$. The proof is based on ideas of Heath-Brown~\cite[Lemma~2]{H79}. 

\begin{proposition}[Increasing property]
\label{prop:increasing property}
For any $J\gg1$, there exists $\theta\in\{1,2\}$ such that 
    \[
    \widetilde{S}(\cX_k,\alpha,N,K) \ll \frac{(JN)^{o(1)}}{J}\widetilde{S}(\cX_k,\alpha,\theta JN,JK). 
    \]
\end{proposition}

Our final recursive bound for $\widetilde{S}$ is implicit in an argument of Heath-Brown~\cite[Theorem~4]{H18} and increases the length of summation. 

\begin{proposition}[Multiplicative property]
\label{prop:multiplicative property}
There exists $\theta\in\{1,2\}$ such that 
    \[
    \widetilde{S}(\cX_k,\alpha,N,K) \ll N^{o(1)}\|\alpha_k\|_{2,N/K}\widetilde{S}(\cX_k,\alpha,\theta N^2,NK)^{1/2}. 
    \]
\end{proposition}

\subsection{Proof of Theorem~\ref{thm:dirichlet} assuming Propositions~\ref{prop:dyadic partition}--\ref{prop:multiplicative property}}

By Propositions~\ref{prop:dyadic partition} and~\ref{prop:linear property}, there exists $N\leq2M$ such that 
    \begin{equation}
    \label{eq:initial decomp}
    S(\cX,\alpha,M,K) \ll (MK)^{o(1)}\sum_{k=0}^\infty\widetilde{S}(\cX_k,\alpha,N,K). 
    \end{equation}
For any $N_0\leq N$, by Proposition~\ref{prop:decreasing property}, there exists $L\leq4N_0$ such that 
    \[
    \widetilde{S}(\cX_k,\alpha,N,K) 
    \ll \frac{NK}{2^{k}}\|\alpha_k\|^2_1+\frac{N}{N_0}\|\alpha_k\|^2_1+\frac{N^{1+o(1)}}{N_0}\widetilde{S}\left(\cX_k,\alpha,L,\frac{N_0K}{N}\right). 
    \]
By Proposition~\ref{prop:multiplicative property}, there exists $\theta_1\in\{1,2\}$ such that 
    \[
    \widetilde{S}\left(\cX_k,\alpha,L,\frac{N_0K}{N}\right) 
    \ll N^{o(1)}\|\alpha_k\|_{2,LN/(N_0K)}\widetilde{S}
    \left(\cX_k,\alpha,\theta_1L^2,\frac{LN_0K}{N}\right)^{1/2}.
    \]
By Proposition~\ref{prop:increasing property} (with $J=N_0^2/L^2$), there exists $\theta_2\in\{1,2\}$ such that 
    \[
    \widetilde{S}\left(\cX_k,\alpha,\theta_1L^2,\frac{LN_0K}{N}\right) \ll \frac{L^{2+o(1)}}{N_0^2}\widetilde{S}\left(\cX_k,\alpha,\theta N_0^2,\frac{N_0^3K}{LN}\right)
    \]
where $\theta=\theta_1^2\theta_2$ satisfies $1\leq\theta\leq8$. Applying Proposition \ref{prop:linear property} (with $J=N_0/L$) to $\widetilde{S}$ on the right-hand side of the above, we obtain 
    \[
    \widetilde{S}\left(\cX_k,\alpha,\theta_1L^2,\frac{LN_0K}{N}\right)
    \ll \frac{L^{1+o(1)}}{N_0}\widetilde{S}
    \left(\cX_k,\alpha,\theta N_0^2,\frac{N_0^2K}{N}\right). 
    \]
Moreover, by Lemma~\ref{lem:l2 norm increasing}, 
    \[
    \|\alpha_k\|_{2,LN/(N_0K)} \ll \left(\frac{N_0}{L}\right)^{1/2}\|\alpha_k\|_{2,N/K}. 
    \]
Combining the above estimates, we have 
    \[
    \widetilde{S}\left(\cX_k,\alpha,L,\frac{N_0K}{N}\right) \ll N^{o(1)}\|\alpha_k\|_{2,N/K}\widetilde{S}\left(\cX_k,\alpha,\theta N_0^2,\frac{N_0^2K}{N}\right)^{1/2}. 
    \]
Therefore, 
    \begin{align*}
    \widetilde{S}(\cX_k,\alpha,N,K) &\ll \frac{NK}{2^{k}}\|\alpha_k\|^2_1+\frac{N}{N_0}\|\alpha_k\|^2_1 \\
    &\quad +\frac{N^{1+o(1)}}{N_0}\|\alpha_k\|_{2,N/K}\widetilde{S}\left(\cX_k,\alpha,\theta N_0^2,\frac{N_0^2K}{N}\right)^{1/2}. 
    \end{align*}

Substituting the above into \eqref{eq:initial decomp} and applying \eqref{eq:norm sums}, we obtain 
    \begin{align*}
    S(\cX,\alpha,M,K) &\ll (MK)^{o(1)}NK\sum_{k=0}^\infty\frac{\|\alpha_k\|_1^2}{2^k}+(MK)^{o(1)}\frac{N}{N_0}\|\alpha\|_1^2 \\
    &\quad +(MK)^{o(1)}\frac{N}{N_0}\sum_{k=0}^\infty\|\alpha_k\|_{2,N/K}\widetilde{S}\left(\cX_k,\alpha,\theta N_0^2,\frac{N_0^2K}{N}\right)^{1/2}. 
    \end{align*}
Moreover, by the Cauchy--Schwarz inequality, 
    \begin{align*}
    S(\cX,\alpha,M,K) &\ll (MK)^{1+o(1)}\sum_{k=0}^\infty\frac{\|\alpha_k\|_1^2}{2^k}+(MK)^{o(1)}\frac{N}{N_0}\|\alpha\|_1^2 \\
    &\quad +(MK)^{o(1)}\frac{N}{N_0}\|\alpha\|_{2,N/K}\left(\sum_{k=0}^\infty\widetilde{S}\left(\cX_k,\alpha,\theta N_0^2,\frac{N_0^2K}{N}\right)\right)^{1/2}. 
    \end{align*}
Since 
    \begin{align*}
    \sum_{k=0}^\infty\widetilde{S}\left(\cX_k,\alpha,\theta N_0^2,\frac{N_0^2K}{N}\right) &\leq \sum_{k=0}^\infty S\left(\cX_k,\alpha,2\theta N_0^2,\frac{N_0^2K}{N}\right) \\ 
    &\leq S\left(\cX,\alpha,2\theta N_0^2,\frac{N_0^2K}{N}\right), 
    \end{align*}
we get
    \begin{align*}
    S(\cX,\alpha,M,K) &\ll (MK)^{1+o(1)}\sum_{k=0}^\infty\frac{\|\alpha_k\|_1^2}{2^k}+(MK)^{o(1)}\frac{N}{N_0}\|\alpha\|_1^2 \\
    &\quad +(MK)^{o(1)}\frac{N}{N_0}\|\alpha\|_{2,N/K}S\left(\cX,\alpha,2\theta N_0^2,\frac{N_0^2K}{N}\right)^{1/2}. 
    \end{align*}

We choose 
    \[
    N_0 = \left(\frac{N}{32}\right)^{1/2}
    \]
so that 
    \begin{align*}
    S(\cX,\alpha,M,K) &\ll (MK)^{1+o(1)}\sum_{k=0}^\infty\frac{\|\alpha_k\|_1^2}{2^k}+(MK)^{o(1)}N^{1/2}\|\alpha\|_1^2 \\
    &\quad +(MK)^{o(1)}N^{1/2}\|\alpha\|_{2,N/K}S\left(\cX,\alpha,\frac{N}{2},\frac{K}{16}\right)^{1/2} \\
    &\leq (MK)^{1+o(1)}\sum_{k=0}^\infty\frac{\|\alpha_k\|_1^2}{2^k}+(MK)^{o(1)}N^{1/2}\|\alpha\|_1^2 \\
    &\quad +(MK)^{o(1)}N^{1/2}\|\alpha\|_{2,N/K}S(\cX,\alpha,M,K)^{1/2}.
    \end{align*}
This implies
    \begin{align*}
    S(\cX,\alpha,M,K) &\ll (MK)^{1+o(1)}\sum_{k=0}^\infty\frac{\|\alpha_k\|_1^2}{2^k}+M^{1/2+o(1)}K^{o(1)}\|\alpha\|_1^2 \\
    &\quad +(MK)^{o(1)}N\|\alpha\|_{2,N/K}^2. 
    \end{align*}
By Lemma~\ref{lem:l2 norm increasing}, 
    \begin{align}
    \label{eq:S bound 0}
    S(\cX,\alpha,M,K) &\ll (MK)^{1+o(1)}\sum_{k=0}^\infty\frac{\|\alpha_k\|_1^2}{2^k}+M^{1/2+o(1)}K^{o(1)}\|\alpha\|_1^2 \\
    &\quad +M^{1+o(1)}K^{o(1)}\|\alpha\|_{2,M/K}^2. \nonumber
    \end{align}

\begin{remark}
While the bound~\eqref{eq:S bound 0} is insufficient for our application to Theorem~\ref{thm:main}, it serves as a preliminary step in the proof of Theorem~\ref{thm:dirichlet}. Moreover, it may be more effective in settings outside the metric Poissonian problem, particularly with alternative choices of parameters. 
\end{remark}

Returning to~\eqref{eq:initial decomp}, we apply Proposition~\ref{prop:multiplicative property} to each $\widetilde{S}$ on the right-hand side, getting 
    \[
    S(\cX,\alpha,M,K) \ll (MK)^{o(1)}\sum_{k=0}^\infty\|\alpha_k\|_{2,N/K}\widetilde{S}(\cX_k,\alpha,\sigma_kN^2,NK)^{1/2}
    \]
for some $\sigma_k\in\{1,2\}$. By the Cauchy--Schwarz inequality, 
    \begin{align*}
    S(\cX,\alpha,M,K)^2 &\ll (MK)^{o(1)}\|\alpha\|_{2,N/K}^2\sum_{k=0}^\infty\widetilde{S}(\cX_k,\alpha,\sigma_kN^2,NK) \\
    &\leq (MK)^{o(1)}\|\alpha\|_{2,N/K}^2\sum_{k=0}^\infty S(\cX_k,\alpha,4N^2,NK) \\
    &\leq (MK)^{o(1)}\|\alpha\|_{2,N/K}^2S(\cX,\alpha,4N^2,NK). 
    \end{align*}
We apply \eqref{eq:S bound 0} to $S(\cX,\alpha,4N^2,NK)$, getting 
    \begin{align*}
    S(\cX,\alpha,4N^2,NK) &\ll (N^3K)^{1+o(1)}\sum_{k=0}^\infty\frac{\|\alpha_k\|_1^2}{2^k}+N^{1+o(1)}K^{o(1)}\|\alpha\|_1^2 \\
    &\quad +N^{2+o(1)}K^{o(1)}\|\alpha\|_{2,N/K}^2
\end{align*}
and results in 
    \begin{align*}
    S(\cX,\alpha,M,K)^2 &\ll M^{o(1)}N^{3+o(1)}K^{1+o(1)}\|\alpha\|_{2,N/K}^2\sum_{k=0}^\infty\frac{\|\alpha_k\|_1^2}{2^k} \\
    &\quad +M^{o(1)}N^{1+o(1)}K^{o(1)}\|\alpha\|_1^2\|\alpha\|_{2,N/K}^2 \\
    &\quad +M^{o(1)}N^{2+o(1)}K^{o(1)}\|\alpha\|_{2,N/K}^4.
    \end{align*}
Finally, by Lemma~\ref{lem:l2 norm increasing}, 
    \begin{align*}
    S(\cX,\alpha,M,K)^2 &\ll M^{1+o(1)}N^{2+o(1)}K^{2+o(1)}\|\alpha\|_{2,M}^2\sum_{k=0}^\infty\frac{\|\alpha_k\|_1^2}{2^k} \\
    &\quad +M^{1+o(1)}N^{o(1)}K^{1+o(1)}\|\alpha\|_1^2\|\alpha\|_{2,M}^2 \\
    &\quad +M^{2+o(1)}N^{o(1)}K^{2+o(1)}\|\alpha\|_{2,M}^4 \\
    &\ll M^{3+o(1)}K^{2+o(1)}\|\alpha\|_{2,M}^2\sum_{k=0}^\infty\frac{\|\alpha_k\|_1^2}{2^k} \\
    &\quad +M^{1+o(1)}K^{1+o(1)}\|\alpha\|_1^2\|\alpha\|_{2,M}^2 \\
    &\quad +M^{2+o(1)}K^{2+o(1)}\|\alpha\|_{2,M}^4.
    \end{align*}
We complete the proof by noting that 
    \[
    \sum_{k=0}^\infty\frac{\|\alpha_k\|_1^2}{2^k} = \sum_{k=0}^{\infty}\left(\sum_{x\in \cX_k}\frac{\alpha(x)}{2^{k/2}} \right)^2 \leq \sum_{k=0}^{\infty}\left(\sum_{x\in \cX_k}\frac{\alpha(x)}{x^{1/2}}\right)^2 \leq \left(\sum_{x\in \cX}\frac{\alpha(x)}{x^{1/2}}\right)^2, 
    \]
and subject to the validity of Propositions~\ref{prop:dyadic partition}--\ref{prop:multiplicative property}. 

\begin{remark}
Due to the loss of information in the chain of inequalities above, we expect that the condition~\eqref{eq:main growth condition} of our main theorem can be weakened. Moreover, an additional apparent loss comes from the application of the Cauchy-Schwarz inequality before invoking Proposition~\ref{prop:decreasing property}. 
\end{remark}

\section{Combinatorial decomposition}
\label{sec:decomposition}

We present some preliminaries for the proofs of Propositions~\ref{prop:dyadic partition} and \ref{prop:linear property}. 

\begin{lemma}
\label{lem:decomp split} 
For any $X_1,X_2\subseteq\cX$, $N_1,N_2\in\NN$ and $y\in\RR$, 
    \begin{align*}
    &\sum_{\substack{n_i\sim N_i\\x_i\in X_i\\|n_1x_1-n_2x_2-y|\leq K}}\alpha(x_1)\alpha(x_2)\beta(n_1)\beta(n_2) \\
    &\quad\ll \widetilde{S}(X_1,\alpha,|\beta|,N_1,2K)^{1/2}\widetilde{S}(X_2,\alpha,|\beta|,N_2,2K)^{1/2}. 
    \end{align*}
\end{lemma}

\begin{proof}
The left-hand side can be written as 
    \begin{equation}
    \label{equation double sum characteristic}
    S := \sum_{\substack{n_i\sim N_i\\x_i\in X_i}}\alpha(x_1)\alpha(x_2)\beta(n_1)\beta(n_2)\1_{[-K,K]}(n_1x_1-n_2x_2-y)
    \end{equation}
where $\1_{[-K,K]}$ denotes the indicator function of the interval $[-K,K]$.
The function 
    \[
    f(x) := \max\left\{1-\left|\frac{x}{2K}\right|,0\right\} 
    \]
satisfies 
    \begin{equation}
    \label{eq:decomp indicator bound}
    \1_{[-K,K]}(x) \ll f(x) \ll \1_{[-2K,2K]}(x), 
    \end{equation}
and its Fourier transform is given by 
    \[
    \widehat{f}(\xi) = \frac{2K\sin(2\pi K\xi)^2}{(2\pi K\xi)^2}. 
    \]
By Fourier inversion and the fact that $\widehat{f}(\xi)\geq0$ for all $\xi\in\RR$, we have
    \begin{align*}
    S &\ll\sum_{\substack{n_i\sim N_i\\x_i\in X_i}}\alpha(x_1)\alpha(x_2)|\beta(n_1)\beta(n_2)|f(n_1x_1-n_2x_2-y)  \\
    &\leq \int_{-\infty}^\infty \widehat{f}(\xi)\left|\sum_{\substack{n_1\sim N_1\\x_1\in X_1}}\alpha(x_1)|\beta(n_1)|e(\xi n_1x_1)\right| \\
    &\quad\times \left|\sum_{\substack{n_2\sim N_2\\x_2\in X_2}}\alpha(x_2)|\beta(n_2)|e(-\xi n_2x_2)\right|d\xi. 
    \end{align*}
By the Cauchy--Schwarz inequality, 
    \begin{align*}
    S^2 &\ll \left(\int_{-\infty}^\infty\widehat{f}(\xi)\left|\sum_{\substack{n_1\sim N_1\\x_1\in X_1}}\alpha(x_1)|\beta(n_1)|e(\xi n_1x_1)\right|^2\,d\xi\right) \\
    &\quad \times\left(\int_{-\infty}^\infty\widehat{f}(\xi)\left|\sum_{\substack{n_2\sim N_2\\x_2\in X_2}}\alpha(x_2)|\beta(n_2)|e(-\xi n_2x_2)\right|^2\,d\xi\right). 
    \end{align*}
Expanding the square and applying~\eqref{eq:decomp indicator bound}, we obtain 
    \begin{align*}
    &\quad \int_{-\infty}^\infty\hat{f}(\xi)\left|\sum_{\substack{n_1\sim N_1\\x_1\in\cX_1}}\alpha(x_1)|\beta(n_1)|e(\xi n_1x_1)\right|^2\,d\xi \\
    &= \sum_{\substack{n_1,n_1'\sim N_1\\x_1,x_1'\in X_1}}\alpha(x_1)\alpha(x_1')|\beta(n_1)\beta(n_1')|\int_{-\infty}^\infty\hat{f}(\xi)e(\xi(n_1x_1-n_1'x_1'))\,d\xi \\
    &= \sum_{\substack{n_1,n_1'\sim N_1\\x_1,x_1'\in\cX_1}}\alpha(x_1)\alpha(x_1')|\beta(n_1)\beta(n_1')|f(n_1x_1-n_1'x_1') \\
    &\ll \sum_{\substack{n_1,n_1'\sim N_1\\x_1,x_1'\in\cX_1}}\alpha(x_1)\alpha(x_1')|\beta(n_1)\beta(n_1')|\1_{[-2K,2K]}(n_1x_1-n_1'x_1') \\
    &= \widetilde{S}(X_1,\alpha,|\beta|,N_1,2K).
    \end{align*}
A similar argument applies to the second integral and we complete the proof. 
\end{proof}

\begin{lemma}
\label{lem:decomp dyadic}
For any $X\subseteq\cX$, there exists $N\le 2M$ such that 
    \[
    S(X,\alpha,M,K) \ll (\log M)^2\widetilde{S}(X,\alpha,N,2K).
    \]
\end{lemma}

\begin{proof}
Splitting $M$ into dyadic intervals and applying Lemma~\ref{lem:decomp split}, we have 
    \begin{align*}
    S(X,\alpha,M,K) &= \sum_{\substack{j,k \\ 2^{j},2^{k}\le M}}\sum_{\substack{n_1\sim2^j\\n_2\sim2^k\\x_1,x_2\in X\\|n_1x_1-n_2x_2|\leq K}}\alpha(x_1)\alpha(x_2) \\
    &\ll \sum_{\substack{j,k \\ 2^{j},2^{k}\le M}}\widetilde{S}(X,\alpha,2^j,2K)^{1/2}\widetilde{S}(X,\alpha,2^k,2K)^{1/2} \\
    &\ll \left(\sum_{\substack{j \\ 2^{j}\le M}}\widetilde{S}(X,\alpha,2^j,2K)^{1/2}\right)^2.
    \end{align*}
The result follows by taking the maximum over $j$ in the above. 
\end{proof}

\begin{lemma}
\label{lem:decomp X dyadic}
There holds 
    \[
    \widetilde{S}(\cX,\alpha,N,K) \ll (1+\log K)\sum_{k=0}^\infty\widetilde{S}(\cX_k,\alpha,N,2K).
    \]
\end{lemma}

\begin{proof}
Recall from \eqref{eq:dyadic S def} that 
    \[
    \widetilde{S}(\cX,\alpha,N,K) = \sum_{\substack{n_1,n_2\sim N\\x_1,x_2\in\cX\\|n_1x_1-n_2x_2|\leq K}}\alpha(x_1)\alpha(x_2). 
    \]
If $x_1,x_2,n_1,n_2$ satisfy the above conditions of summation, then 
    \[
    \left|x_1-\frac{n_2x_2}{n_1}\right| \leq K,
    \]
and hence 
    \[
    \frac{x_2}{2}-K \leq x_1 \leq 2x_2+K.
    \]
This implies that if $x_1\in \cX_k$ and $x_2\in \cX_{\ell}$, then 
    \[
    |k-\ell|\ll 1+\log{K}
    \]
and so 
    \[
    \widetilde{S}(\cX,\alpha,N,K) \leq \sum_{\substack{k,\ell\geq0\\|k-j|\ll1+\log{K}}}\sum_{\substack{n_1,n_2\sim N\\x_1\in\cX_k\\x_2\in\cX_\ell\\|n_1x_1-n_2x_2|\leq K}}\alpha(x_1)\alpha(x_2). 
    \]
Applying Lemma~\ref{lem:decomp split}, we conclude that 
    \begin{align*}
    \widetilde{S}(\cX,\alpha,N,K) &\ll \sum_{\substack{k,\ell\geq0\\|k-j|\ll1+\log{K}}}\widetilde{S}(\cX_k,\alpha,N,2K)^{1/2}\widetilde{S}(\cX_\ell,\alpha,N,2K)^{1/2} \\
    &\ll \sum_{\substack{k,\ell\geq0\\|k-j|\ll1+\log{K}}}\left(\widetilde{S}(\cX_k,\alpha,N,2K)+\widetilde{S}(\cX_\ell,\alpha,N,2K)\right) \\
    &\ll (1+\log{K})\sum_{k=0}^\infty\widetilde{S}(\cX_k,\alpha,N,2K). 
    \qedhere
\end{align*}
\end{proof}

\subsection{Proof of Proposition~\ref{prop:dyadic partition}}

We first apply Lemma~\ref{lem:decomp dyadic}, getting 
    \[
    S(\cX,\alpha,M,K) \ll (\log M)^2\widetilde{S}(\cX,\alpha,N,2K)
    \]
for some $N\leq 2M$. Lemma~\ref{lem:decomp X dyadic} then yields 
    \begin{align*}
    S(\cX,\alpha,M,K) &\ll (\log M)^2(1+\log K)\sum_{k=0}^\infty\widetilde{S}(\cX_k,\alpha,N,4K), 
    \end{align*}
which completes the proof. 

\subsection{Proof of Proposition~\ref{prop:linear property}}

We establish the desired result for arbitrary $X\subseteq\cX$. If $J\leq1$, then trivially we have $\widetilde{S}(X,\alpha,\beta,N,JK)\leq\widetilde{S}(X,\alpha,|\beta|,N,K)$. Hence, we may assume that $J\geq1$. Let 
    \[
    [-JK,JK] = \bigcup_{j=1}^{\lceil4J\rceil}I_j
    \]
be a partition of $[-K,K]$ into subintervals $I_j$ of length at most $K/2$. Then 
    \begin{equation}
    \label{eq:decomp sum partition}
    \widetilde{S}(X,\alpha,\beta,N,JK) \ll \sum_{j=1}^{\lceil4J\rceil}\sum_{\substack{n_1,n_2\sim N\\x_1,x_2\in X\\ n_1x_1-n_2x_2\in I_j}}\alpha(x_1)\alpha(x_2)|\beta(n_1)\beta(n_2)|.
    \end{equation}
If $n_1x_1-n_2x_2\in I_j$, then there exists $y_j$ such that 
    \[
    |n_1x_1-n_2x_2-y_j| \leq \frac{K}{2}. 
    \]
By Lemma~\ref{lem:decomp split}, 
    \begin{align*}
    \sum_{\substack{n_1,n_2\sim N\\x_1,x_2\in X\\n_1x_1-n_2x_2\in I_j}}\alpha(x_1)\alpha(x_2)|\beta(n_1)\beta(n_2)| \ll \widetilde{S}(X,\alpha,|\beta|,N,K), 
    \end{align*}
and the result follows from~\eqref{eq:decomp sum partition}. 

\section{Geometry of numbers}
\label{sec:geometry}

We recall some preliminaries from the geometry of numbers. Given a lattice $\Lambda$ of full rank and a symmetric convex body $B$ in $\RR^n$, the $j$-th successive minimum of $B$ with respect to $\Lambda$, where $1\leq j\leq n$, is defined as 
    \[
    \lambda_j = \inf\{\lambda>0:\lambda B\cap\Lambda\text{ contains $j$ linearly independent points}\}. 
    \]
We denote by $\mu$ the Lebesgue measure. The following is Minkowski's second theorem; see, for example,~\cite[Theorem~3.30]{TV10}. 

\begin{lemma}
\label{lem:Minkowski}
Let $\Lambda$ be a lattice of full rank, $B$ be a symmetric convex body in $\RR^n$, and $\lambda_1,\dots,\lambda_n$ be the successive minima of $B$ with respect to $\Lambda$. Then 
    \[
    \frac{\mu(B)}{\mu(\RR^n/\Lambda)} \ll \frac{1}{\lambda_1\cdots\lambda_n} \ll \frac{\mu(B)}{\mu(\RR^n/\Lambda)}. 
    \]
\end{lemma}

We may bound the number of lattice points $|B\cap\Lambda|$ in terms of the successive minima; see, for example,~\cite[Exercise~3.5.6]{TV10}. 

\begin{lemma}
\label{lem:lattice}
Let $\Lambda$ be a lattice of full rank, $B$ be a symmetric convex body in $\RR^n$, and $\lambda_1,\dots,\lambda_n$ be the successive minima of $B$ with respect to $\Lambda$. Then 
    \begin{align*}
    \prod_{j=1}^n\max\left\{1,\frac{1}{\lambda_j}\right\} \ll |B\cap\Lambda| \ll \prod_{j=1}^n\max\left\{1,\frac{1}{\lambda_j}\right\}. 
    \end{align*}
\end{lemma}

Our main application of the above results is to convex bodies defined by linear inequalities. We now compute the volumes of the relevant bodies. 

\begin{lemma}
\label{lem:volume}
Let $k$ be a non-negative integer, $N,K>0$ and $x_1,x_2\sim 2^{k}$ be real numbers. Then 
    \[
    \mu\left(\{(y_1,y_2)\in[-N,N]^2:|x_1y_1-x_2y_2|\leq K\}\right) \ll \frac{NK}{2^{k}}. 
    \]
\end{lemma}

\begin{proof}
The region is contained in the parallelogram enclosed by the lines 
    \[
    y_2 = \frac{x_1y_1+K}{x_2}, \quad y_2 = \frac{x_1y_1-K}{x_2}, \quad y_2 = N \quad \text{and} \quad y_2 = -N. 
    \]
In particular, $y_2=N$ intersects $y_2=(x_1y_1+K)/x_2$ and $y_2=(x_1y_1-K)/x_2$, respectively, at the points 
    \[
    \left(\frac{Nx_2-K}{x_1},N\right) \quad \text{and} \quad \left(\frac{Nx_2+K}{x_1},N\right). 
    \]
Hence, the volume is bounded by 
    \[
    2N\left(\frac{Nx_2+K}{x_1}-\frac{Nx_2-K}{x_1}\right) = \frac{4NK}{x_1} \ll \frac{NK}{2^k}. \qedhere
    \]
\end{proof}

We also require the following lemma on how the successive minima scale with respect to a convex body. 

\begin{lemma}
\label{lem:scaling}
Let $r>0$, $\Lambda$ a lattice of full rank, $B$ be a symmetric convex body in $\RR^2$, and $\lambda_1,\lambda_2$ be the successive minima of $B$ with respect to $\Lambda$. Then, the successive minima of the convex body 
    \[
    \{rb:b\in B\}
    \]
with respect to $\Lambda$ are given by 
    \[
    \frac{\lambda_1}{r} \quad \text{and} \quad \frac{\lambda_2}{r}. 
    \]
\end{lemma}

\begin{proof}
Let $\lambda_1'$ and $\lambda_2'$ denote the successive minima of $rB$ with respect to $\Lambda$. By symmetry, it's sufficient to show that 
    \[
    \lambda_j' \leq \frac{\lambda_j}{r} \quad \text{for} \quad j=1,2. 
    \]
This follows since $v_1\in\lambda_1B$ and $v_2\in\lambda_2B$ implies that 
    \[
    v_j \in \frac{\lambda_j}{r}(rB). \qedhere
    \]
\end{proof}

\subsection{Proof of Proposition~\ref{prop:decreasing property}}

For $x_1,x_2\in\cX$, define the convex body 
    \[
    B(x_1,x_2,N,K) = \{(y_1,y_2)\in[-2N,2N]^2:|x_1y_1-x_2y_2|\leq K\}. 
    \]
Then 
    \[
    \widetilde{S}(\cX_k,\alpha,N,K) 
    \leq \sum_{x_1,x_2\in\cX_k}\alpha(x_1)\alpha(x_2)
    |B(x_1,x_2,N,K)\cap\ZZ^2|. 
    \]
Let $\lambda_1(x_1,x_2)$ and $\lambda_2(x_1,x_2)$ (which also depend on $N$ and $K$) be the successive minima of $B(x_1,x_2,N,K)$ with respect to $\ZZ^2$. We further split the above summation into 
    \begin{equation}
    \label{eq:geometry S1+S2}
    \widetilde{S}(\cX_k,\alpha,N,K) \leq S_1+S_2+\|\alpha_k\|^2_1
    \end{equation}
where 
    \[
    S_1 := \sum_{\substack{x_1,x_2\in\cX_k\\\lambda_2(x_1,x_2)\leq1}}\alpha(x_1)\alpha(x_2)|B(x_1,x_2,N,K)\cap\ZZ^2| 
    \]
and 
    \[
    S_2 := \sum_{\substack{x_1,x_2\in\cX_k\\\lambda_2(x_1,x_2)>1\\\lambda_1(x_1,x_2)\leq1}}\alpha(x_1)\alpha(x_2)|B(x_1,x_2,N,K)\cap\ZZ^2|, 
    \]
noting that the convex bodies with $\lambda_1(x_1,x_2)>1$ contain exactly one point $(0,0)$ which contribute at most $\|\alpha_k\|_1^2$ to our bound. For $x_1,x_2\in\cX_k$ such that $\lambda_2(x_1,x_2)\leq1$, Lemmata \ref{lem:Minkowski}, \ref{lem:lattice} and \ref{lem:volume} yield 
    \[
    |B(x_1,x_2,N,K)\cap\ZZ^2| \ll \frac{1}{\lambda_1(x_1,x_2)\lambda_2(x_1,x_2)} \ll \frac{KN}{2^{k}}.
    \]
This implies that 
    \begin{equation}
    \label{eq:geometry S1 bound}
    S_1 \ll \frac{KN}{2^{k}}\sum_{\substack{x_1,x_2\in\cX_k\\\lambda_2(x_1,x_2)\leq1}}\alpha(x_1)\alpha(x_2) \leq \frac{KN}{2^k}\|\alpha_k\|_1^2. 
    \end{equation}

Consider next $S_2$. An application of Lemma~\ref{lem:lattice} results in
    \[
    S_2 \ll \sum_{\substack{x_1,x_2\in\cX_k\\\lambda_2(x_1,x_2)>1\\\lambda_1(x_1,x_2)\leq1}}\frac{\alpha(x_1)\alpha(x_2)}{\lambda_1(x_1,x_2)}. 
    \]
Let $r\in(0,1]$ and consider 
    \[
    S_2 \ll \frac{1}{r}\sum_{\substack{x_1,x_2\in\cX_k\\\lambda_2(x_1,x_2)>1\\\lambda_1(x_1,x_2)\leq1}}\alpha(x_1)\alpha(x_2)\frac{r}{\lambda_1(x_1,x_2)}. 
    \]
By Lemma~\ref{lem:scaling}, the successive minima of the convex body 
\begin{align*}
B(x_1,x_2,rN,rK)
\end{align*}
are given by 
    \[
    \frac{\lambda_1(x_1,x_2)}{r} \quad \text{and} \quad \frac{\lambda_2(x_1,x_2)}{r}. 
    \]
Note that 
    \[
    \frac{r}{\lambda_2(x_1,x_2)} < 1. 
    \]
Hence, by Lemma~\ref{lem:lattice}, 
    \[
    \max\left\{1,\frac{r}{\lambda_1(x_1,x_2)}\right\} \ll |B(x_1,x_2,rN,rK)\cap\ZZ^2|, 
    \]
which implies that either 
    \[
    \frac{r}{\lambda_1(x_1,x_2)} \ll 1 \quad \text{or} \quad \frac{r}{\lambda_1(x_1,x_2)} \ll |B(x_1,x_2,rN,rK)\cap\ZZ^2|. 
    \]
It follows that 
    \[
    S_2 \ll \frac{1}{r}\sum_{\substack{x_1,x_2\in\cX_k\\\lambda_2(x_1,x_2)>1\\\lambda_1(x_1,x_2)\leq1}}\alpha(x_1)\alpha(x_2)(1+|B(x_1,x_2,rN,rK)\cap\ZZ^2|).
    \]
Applying the above with 
    \[
    r = \frac{N_0}{N} 
    \]
results in
    \begin{align*}
    S_2 &\ll \frac{N}{N_0}\sum_{\substack{x_1,x_2\in\cX_k\\\lambda_2(x_1,x_2)>1\\\lambda_1(x_1,x_2)\leq1}}\alpha(x_1)\alpha(x_2)\left(1+\left|B\left(x_1,x_2,N_0,\frac{N_0K}{N}\right)\cap\ZZ^2\right|\right) \\
    &\ll \frac{N}{N_0}\|\alpha_k\|_1^2
    +\frac{N}{N_0}\sum_{x_1,x_2\in\cX_k}\alpha(x_1)\alpha(x_2)
    \left|B\left(x_1,x_2,N_0,\frac{N_0K}{N}\right)\cap\ZZ^2\right| \\
    &\ll \frac{N}{N_0}\|\alpha_k\|_1^2
    +\frac{N}{N_0}\sum_{\substack{-2N_0\leq n_1,n_2\leq 2N_0\\x_1,x_2\in\cX_k\\|n_1x_1-n_2x_2|\leq N_0K/N}}\alpha(x_1)\alpha(x_2). 
    \end{align*}
In particular, we have 
    \begin{align*}
    \sum_{\substack{-2N_0\leq n_1,n_2\leq 2N_0\\x_1,x_2\in\cX_k\\|n_1x_1-n_2x_2|\leq N_0K/N}}\alpha(x_1)\alpha(x_2) &\ll \sum_{x_1,x_2\in\cX_k}\alpha(x_1)\alpha(x_2) \\
    &\quad +\sum_{\substack{1\leq n_1,n_2\leq 2N_0\\x_1,x_2\in\cX_k\\|n_1x_1-n_2x_2|\leq N_0K/N}}\alpha(x_1)\alpha(x_2) \\
    &\quad +\sum_{\substack{1\leq n_1,n_2\leq 2N_0\\x_1,x_2\in\cX_k\\n_1x_1+n_2x_2\leq N_0K/N}}\alpha(x_1)\alpha(x_2) \\
    &\ll \|\alpha_k\|_1^2+S\left(\cX_k,\alpha,2N_0,\frac{N_0K}{N}\right). 
    \end{align*}
Moreover, by Lemma~\ref{lem:decomp dyadic}, there exists $L\le 4N_0$ such that 
    \[
    S\left(\cX_k,\alpha,2N_0,\frac{N_0K}{N}\right) \ll (\log N_0)^2\widetilde{S}\left(\cX_k,\alpha,L,\frac{2N_0K}{N}\right). 
    \]
Therefore, 
    \[
    S_2 \ll \frac{N}{N_0}\|\alpha_k\|_1^2+\frac{N(\log N_0)^2}{N_0}S\left(\cX_k,\alpha,L,\frac{2N_0K}{N}\right). 
    \]
Applying Proposition~\ref{prop:linear property} to the above, and combining the resulting estimate with~\eqref{eq:geometry S1+S2} and \eqref{eq:geometry S1 bound}, we complete the proof. 

\section{Mean values of Dirichlet polynomials}
\label{sec:dirichlet}

The key observation underlying our approach is that the lattice point counting problem~\eqref{eq:dyadic S def} can be written, up to an absolute constant, as the mean value of a Dirichlet polynomial. The following lemma is a variant of~\cite[Lemma~2.1]{W89}. 

\begin{lemma}
\label{lem:sum integral} 
For any $T\gg1$ and $\beta(n)\ll1$, 
    \begin{align*}
    &\widetilde{S}\left(\cX_k,\alpha,\beta,N,\frac{2^kN}{T}\right) \\
    &\qquad \ll \frac{1}{T}\int_{-T}^T\left|\sum_{x\in \cX_k}\frac{\alpha(x)}{x^{it}}\sum_{n\sim N}\frac{\beta(n)}{n^{it}}\right|^2dt \ll \widetilde{S}\left(\cX_k,\alpha,|\beta|,N,\frac{2^kN}{T}\right). 
    \end{align*}
\end{lemma}

\begin{proof}
Let $K=2^kN/T$. It is sufficient to show that there exists absolute constants $c$ and $C$ such that 
    \begin{align*}
    &\widetilde{S}(\cX_k,\alpha,|\beta|,N,cK) \\
    &\qquad \ll \frac{1}{T}\int_{-T}^T\left|\sum_{x\in \cX_k}\frac{\alpha(x)}{x^{it}}\sum_{n\sim N}\frac{\beta(n)}{n^{it}}\right|^2dt 
    \ll \widetilde{S}(\cX_k,\alpha,|\beta|,N,CK), 
    \end{align*}
since the result then follows from Proposition~\ref{prop:linear property} upon noting that 
    \[
    \widetilde{S}(\cX_k,\alpha,\beta,N,K) \ll \frac{1}{c}\widetilde{S}\left(\cX_k,\alpha,|\beta|,N,cK\right) \ll \widetilde{S}\left(\cX_k,\alpha,|\beta|,N,cK\right)
    \]
and 
    \[
    \widetilde{S}(\cX_k,\alpha,|\beta|,N,CK) \ll \widetilde{S}(\cX_k,\alpha,|\beta|,N,K). 
    \]
Consider the function 
    \begin{equation}
    \label{eq:w def}
    w(x) := \max\{1-|x|,0\}, 
    \end{equation}
which has Fourier transform 
    \[
    \widehat{w}(\xi) = \left(\frac{\sin(\pi\xi)}{\pi\xi}\right)^2. 
    \]
We have 
    \[
    \widehat{w}\left(\frac{t}{2T}\right) 
     \gg 1 \quad \text{for} \quad |t|\le T
    \]
so that 
    \begin{align*}
    \int_{-T}^{T}\left|\sum_{x\in \cX_k}\frac{\alpha(x)}{x^{it}}\sum_{n\sim N}\frac{\beta(n)}{n^{it}}\right|^2dt 
    &\ll \int_{-\infty}^\infty\widehat{w}\left(\frac{t}{2T}\right)\left|\sum_{x\in \cX_k}\frac{\alpha(x)}{x^{it}}\sum_{n\sim N}\frac{\beta(n)}{n^{it}}\right|^2dt.
    \end{align*}
Expanding the square, rearranging and using Fourier inversion, we get
    \begin{align*}
    &\quad \int_{-T}^{T}\left|\sum_{x\in \cX_k}\frac{\alpha(x)}{x^{it}}\sum_{n\sim N}\frac{\beta(n)}{n^{it}}\right|^2dt \\
    &\ll T\sum_{x_1,x_2\in \cX_k}\sum_{n_1,n_2\sim N}\alpha(x_1)\alpha(x_2)\beta(n_1)\overline{\beta(n_2)}w\left(\frac{T}{\pi}\log\left(\frac{n_1x_1}{n_2x_2}\right)\right) \\
    &\ll T\sum_{\substack{n_1,n_2\sim N\\x_1,x_2\in\cX_k\\|\log((n_1x_1)/(n_2x_2))|\ll 1/T}}\alpha(x_1)\alpha(x_2)|\beta(n_1)\beta(n_2)|. 
    \end{align*}
If $n_i,x_i$ satisfy the above condition of summation, then 
    \[
    n_1x_1-n_2x_2 \ll \frac{2^{k}N}{T} = K. 
    \]
We conclude that 
    \begin{align*}
    \int_{-T}^{T}\left|\sum_{x\in \cX_k}\frac{\alpha(x)}{x^{it}}\sum_{n\sim N}\frac{\beta(n)}{n^{it}}\right|^2dt 
    &\ll T\widetilde{S}(\cX_k,\alpha,|\beta|,N,CK) 
    \end{align*}
for some absolute constant $C$. 

On the other hand, if $|n_1x_1-n_2x_2|\leq 2^kN/T$, then 
    \[
    \left|\frac{n_1x_1-n_2x_2}{n_2x_2}\right| \ll \frac{1}{T}, 
    \]
and so 
    \[
    \left|\log\left(\frac{n_1x_1}{n_2x_2}\right)\right| \ll \frac{1}{T}. 
    \]
Thus, for some absolute constant $c$, 
    \begin{align*}
    &\quad \widetilde{S}(\cX_k,\alpha,|\beta|,N,cK) \\
    &\leq \sum_{\substack{n_1,n_2\sim N\\x_1,x_2\in\cX_k\\|\log(n_1x_1/(n_2x_2))|\leq 2\pi/T}}\alpha(x_1)\alpha(x_2)|\beta(n_1)\beta(n_2)| \\
    &\ll  \sum_{x_1,x_2\in\cX_k}\sum_{n_1,n_2\sim N}\alpha(x_1)\alpha(x_2)|\beta(n_1)\beta(n_2)|\widehat{w}\left(\frac{T}{2\pi}\log\left(\frac{n_1x_1}{n_2x_2}\right)\right).
    \end{align*}
By Fourier inversion and a change of variable in the resulting integral, we get
    \begin{align*}
    &\quad \widetilde{S}(\cX_k,\alpha,|\beta|,N,cK) \\
    &\ll \frac{1}{T}\sum_{x_1,x_2\in\cX_k}\sum_{n_1,n_2\sim N}\alpha(x_1)\alpha(x_2)|\beta(n_1)\beta(n_2)| \\
    &\quad \times\int_{-\infty}^\infty w\left(\frac{t}{T}\right)e\left(\frac{t}{2\pi}\log\left(\frac{n_1x_1}{n_2x_2}\right)\right)dt \\
    &\ll \frac{1}{T}\int_{-\infty}^\infty w\left(\frac{t}{T}\right)\sum_{x_1,x_2\in\cX_k}\sum_{n_1,n_2\sim N}\frac{\alpha(x_1)\alpha(x_2)|\beta(n_1)\beta(n_2)|}{(n_1x_1)^{it}(n_2x_2)^{-it}}dt.
\end{align*}
After rearranging and using~\eqref{eq:w def}, we arrive at 
    \[
    \widetilde{S}(\cX_k,\alpha,\beta,N,cK) \ll \frac{1}{T}\int_{-T}^{T}\left|\sum_{x\in\cX_k}\frac{\alpha(x)}{x^{it}}\sum_{n\sim N}\frac{\beta(n)}{n^{it}}\right|^2dt. \qedhere
    \]
\end{proof}

In the proof of the following lemma, we adopt an idea of Heath-Brown from~\cite[Lemma~2]{H79}. 

\begin{lemma}
\label{lem:increase dyadic range}
For any $J\gg1$, there exists $\theta\in\{1,2\}$ such that 
    \[
    \int_{-T}^T\left|\sum_{x\in\cX_k}\frac{\alpha(x)}{x^{it}}
    \sum_{n\sim N}\frac{1}{n^{it}}\right|^2dt 
    \ll \frac{(JN)^{o(1)}}{J}\int_{-T}^T
    \left|\sum_{x\in\cX_k}\frac{\alpha(x)}{x^{it}}
    \sum_{m\sim\theta JN}\frac{1}{m^{it}}\right|^2dt. 
    \]
\end{lemma}
\begin{proof}
Let $q$ be a prime satisfying 
    \[
    4J\le q \le 8J.
    \]
For a character $\chi$ modulo $q$, we define 
    \begin{align*}
    \sigma(\chi) &= \int_{-T}^T
    \left|\sum_{x\in\cX_k}\frac{\alpha(x)}{x^{it}}
    \sum_{n\sim N}\frac{1}{n^{it}}\sum_{j\sim J}
    \frac{\chi(j)}{j^{it}}\right|^2dt \\
    &= \int_{-T}^T\left|\sum_{x\in\cX_k}
    \frac{\alpha(x)}{x^{it}}\sum_{JN<m\leq4JN}
    \frac{a_\chi(m)}{m^{it}}\right|^2dt 
    \end{align*}
where 
    \[
    a_\chi(m) := \sum_{\substack{n\sim N\\j\sim J\\m=jn}}\chi(j). 
    \]
We have
    \[
   |a_{\chi}(m)|\le \max_{JN<m\leq4JN}
    \sum_{\substack{n\sim N\\j\sim J\\nj=m}}1 \ll (JN)^{o(1)}
    \]
so that
    \[
    \sigma(\chi) \ll N^{o(1)}\int_{-T}^T\left|
    \sum_{x\in\cX_k}\frac{\alpha(x)}{x^{it}}
    \sum_{JN<m\leq4JN}\frac{\beta_\chi(m)}{m^{it}}\right|^2dt
    \]
for some $\beta_\chi(m)\ll1$. Moreover, we split the second summation into dyadic ranges, getting 
    \begin{align*}
    \sigma(\chi) &\ll N^{o(1)}\int_{-T}^T\left|\sum_{x\in\cX_k}\frac{\alpha(x)}{x^{it}}
    \left(\sum_{m\sim JN}\frac{\beta_\chi(m)}{m^{it}}
    +\sum_{m\sim2JN}\frac{\beta_\chi(m)}{m^{it}}\right)\right|^2dt \\
    &\ll N^{o(1)}\int_{-T}^T\left|\sum_{x\in\cX_k}
    \frac{\alpha(x)}{x^{it}}
    \sum_{m\sim\theta JN}\frac{\beta_\chi(m)}{m^{it}}\right|^2dt
    \end{align*}
for some $\theta\in\{1,2\}$. By two successive applications of Lemma~\ref{lem:sum integral}, 
    \begin{align*}
    \sigma(\chi) &\ll N^{o(1)}\int_{-T}^T\left|
    \sum_{x\in\cX_k}\frac{\alpha(x)}{x^{it}}
    \sum_{m\sim\theta JN}\frac{1}{m^{it}}\right|^2dt. 
    \end{align*}
We note that the right-hand side is now independent of $\chi$. Summing over all characters $\chi$ modulo $q$, we obtain 
    \begin{equation} \label{eq:sum chi upper bound} 
    \sum_\chi\sigma(\chi) \leq \varphi(q)N^{o(1)}\int_{-T}^T\left|\sum_{x\in\cX_k}\frac{\alpha(x)}{x^{it}}\sum_{m\sim JN}\frac{1}{m^{it}}\right|^2dt. 
    \end{equation}

On the other hand, since $q>J$, we have 
    \[
    \sum_\chi\left|\sum_{j\sim J}\frac{\chi(j)}{j^{it}}\right|^2 
    \geq \varphi(q)J 
    \]
and hence 
    \[
    \sum_\chi\sigma(\chi) \geq \varphi(q)J\int_{-T}^T\left|\sum_{x\in\cX_k}\frac{\alpha(x)}{x^{it}}\sum_{n\sim N}\frac{1}{n^{it}}\right|^2dt. 
    \]
Combining this with \eqref{eq:sum chi upper bound}, the result follows. 
\end{proof}

The following lemma can be established in a manner similar to the proof of Lemma~\ref{lem:sum integral} and is applied in a more general context. 

\begin{lemma}
\label{lem:mean value theorem}
For any $\alpha_m\in\CC$ and $x_m\in\RR$ where $m\leq M$, 
    \begin{align*}
    &\sum_{\substack{m_1,m_2\leq M\\2T|x_{m_1}-x_{m_2}|\leq1}}|\alpha_{m_1}\alpha_{m_2}| \\
    &\qquad \ll \frac{1}{T}\int_{-T}^T\left|\sum_{m\leq M}\alpha_me(x_mt)\right|^2\,dt \ll \sum_{\substack{m_1,m_2\leq M\\2T|x_{m_1}-x_{m_2}|\leq1}}|\alpha_{m_1}\alpha_{m_2}|. 
    \end{align*}
\end{lemma}

\subsection{Proof of Proposition~\ref{prop:increasing property}}

By Lemma \ref{lem:sum integral}, 
    \[
    \widetilde{S}(\cX_k,\alpha,N,K) 
    \ll \frac{K}{2^kN}\int_{-2^kN/K}^{2^kN/K}\left|
    \sum_{x\in\cX_k}\frac{\alpha(x)}{x^{it}}
    \sum_{n\sim N}\frac{1}{n^{it}}\right|^2dt. 
    \]
By Lemma \ref{lem:increase dyadic range}, there exists $\theta\in\{1,2\}$ such that 
    \[
    \widetilde{S}(\cX_k,\alpha,N,K) 
    \ll \frac{(JN)^{o(1)}K}{2^kJN}\int_{-2^kN/K}^{2^kN/K}
    \left|\sum_{x\in\cX_k}\frac{\alpha(x)}{x^{it}}
    \sum_{n\sim\theta JN}\frac{1}{n^{it}}\right|^2dt. 
    \]
Applying Proposition~\ref{prop:linear property} and Lemma \ref{lem:sum integral}, we obtain 
    \[
    \widetilde{S}(\cX_k,\alpha,N,K) \ll \frac{(JN)^{o(1)}}{J}
    \widetilde{S}(\cX_k,\alpha,\theta JN,JK). 
    \]
This completes the proof.

\subsection{Proof of Proposition~\ref{prop:multiplicative property}}

By Lemma~\ref{lem:sum integral} and the Cauchy--Schwarz inequality, 
    \begin{align*}
    \frac{2^kN}{K}\widetilde{S}(\cX_k,\alpha,N,K) &\ll \left(\int_{-2^kN/K}^{2^kN/K}\left|\sum_{x\in\cX_k}\frac{\alpha(x)}{x^{it}}\right|^2dt\right)^{1/2} \\
    &\quad \times\left(\int_{-2^kN/K}^{2^kN/K}\left|\sum_{x\in\cX_k}\frac{\alpha(x)}{x^{it}}\left(\sum_{n\sim N}\frac{1}{n^{it}}\right)^2\right|^2dt\right)^{1/2}.
    \end{align*}
Applying Lemma \ref{lem:mean value theorem} and then Lemma~\ref{lem:l2 norm increasing}, we have 
    \begin{align*}
    \frac{K}{2^{k}N}\int_{-2^kN/K}^{2^kN/K}\left|\sum_{x\in\cX_k}
    \frac{\alpha(x)}{x^{it}}\right|^2dt 
    &\ll \sum_{\substack{x_1,x_2\in\cX_k\\
    |\log(x_1/x_2)|\ll  K/(2^kN)}}\alpha(x_1)\alpha(x_2) \\
    &\ll \|\alpha_k\|^2_{2,N/K}.
    \end{align*}
This implies that 
    \begin{align}
    \label{eq:S second integral}
    &\quad \widetilde{S}(\cX_k,\alpha,N,K) \\
    &\ll \|\alpha_k\|_{2,N/K}\left(\frac{K}{2^kN}
    \int_{-2^kN/K}^{2^kN/K}\left|\sum_{x\in\cX_k}
    \frac{\alpha(x)}{x^{it}}\left(\sum_{n\sim N}
    \frac{1}{n^{it}}\right)^2\right|^2dt\right)^{1/2}. \nonumber
    \end{align}
Since 
    \[
    \max_{N^2<m\leq4N^2}\sum_{\substack{a,b\sim N\\ab=m}}1 
    \ll N^{o(1)}, 
    \]
we have 
    \begin{align*}
    &\quad \int_{-2^kN/K}^{2^kN/K}\left|\sum_{x\in\cX_k}
    \frac{\alpha(x)}{x^{it}}\left(\sum_{n\sim N}
    \frac{1}{n^{it}}\right)^2\right|^2dt \\
    &\ll N^{o(1)}\int_{-2^kN/K}^{2^kN/K}\left|\sum_{x\in\cX_k}
    \frac{\alpha(x)}{x^{it}}\sum_{N<m\leq4N^2}
    \frac{\beta(m)}{m^{it}}\right|^2dt 
    \end{align*}
for some $\beta(m)\ll1$. Then, we split the second sum into dyadic ranges and apply Minkowski's inequality to get 
    \begin{align*}
    &\quad \int_{-2^kN/K}^{2^kN/K}\left|\sum_{x\in\cX_k}
    \frac{\alpha(x)}{x^{it}}\left(\sum_{n\sim N}
    \frac{1}{n^{it}}\right)^2\right|^2dt \\
    &\ll N^{o(1)}\int_{-2^kN/K}^{2^kN/K}\left|\sum_{x\in\cX_k}
    \frac{\alpha(x)}{x^{it}}\sum_{m\sim\theta N^2}
    \frac{\beta(m)}{m^{it}}\right|^2dt 
    \end{align*}
for some $\theta\in\{1,2\}$. Applying Lemma \ref{lem:sum integral} and then Proposition~\ref{prop:linear property} gives 
    \begin{align*}
   \frac{K}{2^{k}N} \int_{-2^kN/K}^{2^kN/K}\left|\sum_{x\in\cX_k}
    \frac{\alpha(x)}{x^{it}}\left(\sum_{n\sim N}
    \frac{1}{n^{it}}\right)^2\right|^2dt 
    \ll N^{o(1)}\widetilde{S}
    (\cX_k,\alpha,\theta N^2,NK). 
    \end{align*}
Substituting the above to \eqref{eq:S second integral}, we conclude that 
    \[
    \widetilde{S}(\cX_k,\alpha,N,K) 
    \ll N^{o(1)}\|\alpha_k\|_{2,N/K}
    \widetilde{S}(\cX_k,\alpha,\theta N^2,NK)^{1/2}. 
    \]
This completes the proof. 

\section{The metric Poissonian property}
\label{sec:poissonian}

In this section, we prove Theorem~\ref{thm:main} and Corollaries~\ref{cor:energy only} and~\ref{cor:short interval}. We recall the main reduction to the lattice point counting problem due to Rudnick and Technau, stated as Theorem~\ref{thm:RT}. A sequence $(x_n)_{n\in\NN}$ of real numbers has metric Poissonian pair correlation if there exists $\delta>0$ such that for sufficiently small $\varepsilon>0$, 
    \begin{equation}
    \label{eq:RT condition 1}
    \sum_{\substack{1\leq m\leq N^{1+\varepsilon}\\1\leq n_1\neq n_2\leq N\\m|x_{n_1}-x_{n_2}|<N^\varepsilon}}1 \ll N^{2-\delta} 
    \end{equation}
and 
    \begin{equation}
    \label{eq:RT condition 2}
    \sum_{\substack{1\leq |m_1|,|m_2|\leq N^{1+\varepsilon}\\1\leq n_1\neq n_2\leq N\\1\leq n_3\neq n_4\leq N\\|m_1(x_{n_1}-x_{n_2})-m_2(x_{n_3}-x_{n_4})|<N^\varepsilon}}1 \ll N^{4-\delta}. 
    \end{equation}
To establish \eqref{eq:RT condition 2}, we apply Theorem~\ref{thm:dirichlet} with 
    \begin{equation}
    \label{eq:X def again}
    \cX := \{x_{n_2}-x_{n_1}:1\leq n_1<n_2\leq N\} 
    \end{equation}
and 
    \begin{equation}
    \label{eq:alpha def again}
    \alpha(x) := |\{1\leq n_1<n_2\leq N:x_{n_2}-x_{n_1}=x\}|. 
    \end{equation}
Our next result carries out the details of this reduction, resulting in a condition stated in terms of the sums 
    \begin{equation}
    \label{eq:S def again}
    S(\cX,\alpha,M,K) = \sum_{\substack{1\leq m_1,m_2\leq M\\ x_1,x_2\in \cX \\|m_1x_1-m_2x_2|<K}}\alpha(x_1)\alpha(x_2), 
    \end{equation}
as defined in \eqref{eq:S def}. This is similar to some manipulations occurring in~\cite[Section~5]{AEM21}. 

\begin{proposition}
\label{prop:RT condition bound}
Let $(x_n)_{n\in\NN}$ be a sequence of positive real numbers for which there exist $c>0$ and $\eta\in(0,1]$ such that 
    \begin{equation}
    \label{eq:spacing condition again}
    x_{n+1}-x_n > \frac{c}{n^{1-\eta}}. 
    \end{equation}
Suppose that $\cX$ and $\alpha(x)$ are defined as in~\eqref{eq:X def again} and~\eqref{eq:alpha def again}, respectively. Assume that there exists $\delta>0$ such that for sufficiently small $\varepsilon>0$, 
    \begin{equation}
    \label{eq:S bound assumption}
    S(\cX,\alpha,N^{1+\varepsilon},N^\varepsilon) \ll N^{4-\delta}. 
    \end{equation}
Then $(x_n)_{n\in\NN}$ is metric Poissonian. 
\end{proposition}

\begin{proof}
By Theorem~\ref{thm:RT}, it is sufficient to establish the estimates~\eqref{eq:RT condition 1} and~\eqref{eq:RT condition 2} (with potentially different values of $\delta$). The former follows from \eqref{eq:spacing condition again} since 
    \begin{align*}
    \sum_{\substack{1\leq m\leq N^{1+\varepsilon}\\1\leq n_1\neq n_2\leq N\\m|x_{n_1}-x_{n_2}|<N^\varepsilon}}1 &\ll \sum_{\substack{1\leq m\leq N^{1+\varepsilon}\\1\leq n_1<n_2\leq N\\m(x_{n_2}-x_{n_1})<N^\varepsilon}}1 \\
    &\ll \sum_{\substack{1\leq m\leq N^{1+\varepsilon}\\1\leq n_1<n_2\leq N\\m(n_2-n_1)N^{-1+\eta}<N^\varepsilon}}1 \\
    &\ll N\sum_{\substack{1\leq m\leq N^{1+\varepsilon}\\1\leq n\leq N\\mn<N^{1-\eta+\varepsilon}}}1 \\
    &\ll N^{2-\eta+2\varepsilon}. 
    \end{align*}

It remains to deal with the condition~\eqref{eq:RT condition 2}. By~\eqref{eq:S bound assumption}, it is sufficient to show that 
    \begin{equation}
    \label{eq:RT condition sums}
    \sum_{\substack{1\leq |m_1|,|m_2|\leq N^{1+\varepsilon}\\1\leq n_1\neq n_2\leq N\\1\leq n_3\neq n_4\leq N\\|m_1(x_{n_1}-x_{n_2})-m_2(x_{n_3}-x_{n_4})|<N^\varepsilon}}1 \ll S(\cX,\alpha,N^{1+\varepsilon},N^\varepsilon)+N^{4-\delta}. 
    \end{equation}
By~\eqref{eq:alpha def again} and~\eqref{eq:S def again}, we have 
    \[
    S(\cX,\alpha,N^{1+\varepsilon},N^{\varepsilon}) = \sum_{\substack{1\leq m_1,m_2\leq N^{1+\varepsilon}\\1\leq n_1<n_2\leq N\\1\leq n_3<n_4\leq N\\|m_1(x_{n_2}-x_{n_1})-m_2(x_{n_4}-x_{n_3})|<N^{\varepsilon}}}1, 
    \]
which differs from the left-hand side of~\eqref{eq:RT condition sums} only in conditions of summation. 

By a change of variable $m_j\mapsto-m_j$ if $m_j<0$, and since $(x_n)_{n\in\NN}$ is a strictly increasing sequence, it is readily seen that 
    \begin{align*}
    \sum_{\substack{1\leq |m_1|,|m_2|\leq N^{1+\varepsilon}\\1\leq n_1\neq n_2\leq N\\1\leq n_3\neq n_4\leq N\\|m_1(x_{n_1}-x_{n_2})-m_2(x_{n_3}-x_{n_4})|<N^\varepsilon}}1 &\ll S(\cX,\alpha,N^{1+\varepsilon},N^{\varepsilon}) \\
    &\quad +\sum_{\substack{1\leq m_1,m_2\leq N^{1+\varepsilon}\\1\leq n_1<n_2\leq N\\1\leq n_3<n_4\leq N\\m_1(x_{n_2}-x_{n_1})+m_2(x_{n_4}-x_{n_3})<N^{\varepsilon}}}1. 
    \end{align*}
The assumption~\eqref{eq:spacing condition again} implies that 
    \begin{align*}
    \sum_{\substack{1\leq m_1,m_2\leq N^{1+\varepsilon}\\1\leq n_1<n_2\leq N\\1\leq n_3<n_4\leq N\\m_1(x_{n_2}-x_{n_1})-m_2(x_{n_4}-x_{n_3})<N^{\varepsilon}}}1 &\ll \sum_{\substack{1\leq m_1,m_2\leq N^{1+\varepsilon}\\
    1\leq n_1<n_2\leq N\\1\leq n_3<n_4\leq N\\
    m_1(n_2-n_1)+m_2(n_4-n_3)<N^{1-\eta+\varepsilon}}}1 \\
    &\leq N^{2}\sum_{\substack{1\leq m_1,m_2\leq N^{1+\varepsilon}\\m_1+m_2<N^{1-\eta+\varepsilon}}}1 \\
    &\ll N^{4-2\eta+2\varepsilon}, 
    \end{align*}
which establishes~\eqref{eq:RT condition sums}, after redefining $\delta$ if necessary. 
\end{proof}

\subsection{Proof of Theorem~\ref{thm:main}}

With notation as in Proposition~\ref{prop:RT condition bound}, it is sufficient to show that \eqref{eq:S bound assumption} holds; that is, 
    \[
    S(\cX,\alpha,N^{1+\varepsilon},N^\varepsilon) \ll N^{4-\delta}. 
    \]
We first note that Theorem~\ref{thm:dirichlet} applies only with $\cX\subseteq[1,\infty)$, while the choice~\eqref{eq:X def again} of $\cX$ is potentially contained in $(0,\infty)$. We thus consider partitioning $\cX$ into 
    \[
    \cX_- := \cX\cap(0,1) \quad \text{and} \quad \cX_+ := \cX\cap[1,\infty). 
    \]
By Lemma \ref{lem:decomp split}, 
    \begin{align*}
    &\sum_{\substack{1\leq m_1,m_2\leq N^{1+\varepsilon}\\x_1\in\cX_-\\x_2\in\cX_+\\|m_1x_1-m_2x_2|<N^\varepsilon}}\alpha(x_1)\alpha(x_2) \ll S(\cX_-,\alpha,N^{1+\varepsilon},N^\varepsilon)^{1/2}S(\cX_+,\alpha,N^{1+\varepsilon},N^\varepsilon)^{1/2}. 
    \end{align*}
Thus, \eqref{eq:S bound assumption} follows upon establishing both 
    \begin{equation}
    \label{eq:S bound small}
    S(\cX_-,\alpha,N^{1+\varepsilon},N^\varepsilon) \ll N^{4-\delta}
    \end{equation}
and 
    \begin{equation}
    \label{eq:S bound large}
    S(\cX_+,\alpha,N^{1+\varepsilon},N^\varepsilon) \ll N^{4-\delta}. 
    \end{equation}

We first deal with the former. By the spacing condition \eqref{eq:main spacing condition}, we have 
    \begin{align*}
    S(\cX_-,\alpha,N^{1+\varepsilon},N^\varepsilon) &= \sum_{\substack{1\leq m_1\leq N^{1+\varepsilon}\\x_1,x_2\in\cX\cap(0,1)}}\alpha(x_1)\alpha(x_2)\sum_{\substack{1\leq m_2\leq N^{1+\varepsilon}\\|m_1x_1-m_2x_2|<N^\varepsilon}}1 \\
    &\leq \sum_{\substack{1\leq m_1\leq N^{1+\varepsilon}\\x_1,x_2\in\cX\cap(0,1)}}\alpha(x_1)\alpha(x_2)\sum_{\substack{1\leq m_2\leq N^{1+\varepsilon}\\|m_1x_1/x_2-m_2|<N^{1-\eta+\varepsilon}}}1 \\
    &\ll N^{2-\eta+2\varepsilon}\left(\sum_{x\in\cX\cap(0,1)}\alpha(x)\right)^2.
    \end{align*}
Since
    \[
    \sum_{x\in\cX\cap(0,1)}\alpha(x) = \sum_{\substack{1\leq n_1<n_2\leq N\\x_{n_2}-x_{n_1}<1}}1 \leq \sum_{\substack{1\leq n_1<n_2\leq N\\n_2-n_1<N^{1-\eta}}}1 \ll N^{2-\eta}, 
    \]
we have 
    \[
    S(\cX_-,\alpha,N^{1+\varepsilon},N^\varepsilon) \ll N^{6-3\eta+2\varepsilon}. 
    \]
By the choice $\eta>2/3$, there exists $\delta'>0$ such that $\eta=2/3+\delta'$, and the estimate \eqref{eq:S bound small} follows. 

Consider next~\eqref{eq:S bound large}. Applying Theorem~\ref{thm:dirichlet}, we get 
    \begin{align*}
    S(\cX_+,\alpha,N^{1+\varepsilon},N^\varepsilon) &\ll N^{3/2+o(1)}\|\alpha\|_{2,N}\sum_{x\in\cX_{+}}\frac{\alpha(x)}{x^{1/2}} \\
    &\quad +N^{1/2+o(1)}\|\alpha\|_1\|\alpha\|_{2,N}+N^{1+o(1)}\|\alpha\|_{2,N}^2. 
    \end{align*}
Here, the norms are defined as in \eqref{eq:l1 norm def} and \eqref{eq:l2 norm def}, respectively, with the restriction to $\cX_+$. Note also that $\|\alpha\|_{2,N^{1+\varepsilon}}\leq\|\alpha\|_{2,N}$. It remains to establish estimates for $\|\alpha\|_1$, $\|\alpha\|_{2,N}$ and the weighted sum over $\cX$. By definition \eqref{eq:l1 norm def}, we have 
    \[
    \|\alpha\|_1 = \sum_{x\in\cX}\sum_{\substack{1\leq n_1<n_2\leq N\\x_{n_2}-x_{n_1}=x}}1 \leq N^2. 
    \]
Similarly, by definition \eqref{eq:l2 norm def}, we have 
    \begin{align*}
    \|\alpha\|_{2,N}^2 &\ll \sum_{\substack{1\leq n_1<n_2\leq N\\1\leq n_3<n_4\leq N\\|x_{n_1}-x_{n_2}-x_{n_3}+x_{n_4}|<1/N}}1 \ll E^*_{N,1/N}. 
    \end{align*}
For the weighted sum, we have 
    \[
    \sum_{x\in\cX_{+}}\frac{\alpha(x)}{x^{1/2}} = \sum_{\substack{1\leq n_1<n_2\leq N \\ x_{n_2}-x_{n_1}\ge 1}}\frac{1}{(x_{n_2}-x_{n_1})^{1/2}}. 
    \]
It follows that 
    \begin{align*}
    S(\cX_+,\alpha,N^{1+\varepsilon},N^\varepsilon) &\ll N^{3/2+o(1)}(E^*_{N,1/N})^{1/2}\sum_{\substack{1\leq n_1<n_2\leq N \\ x_{n_2}-x_{n_1}\ge 1}}\frac{1}{(x_{n_2}-x_{n_1})^{1/2}} \\
    &\quad +N^{5/2+o(1)}(E^*_{N,1/N})^{1/2}+N^{1+o(1)}E^*_{N,1/N}. 
    \end{align*}
The estimate \eqref{eq:S bound large} follows from the assumptions \eqref{eq:main energy condition} and \eqref{eq:main growth condition}; that is, 
    \[
    S(\cX_+,\alpha,N^{1+\varepsilon},N^\varepsilon) \ll N^{4-\delta/2+o(1)}+N^{4-\kappa/2-\delta/2+o(1)}+N^{4-\kappa-\delta+o(1)}. 
    \]
This completes the proof. 

\subsection{Improved energy bound: proof of Corollary~\ref{cor:energy only}}

Our goal is to bound the left-hand side of~\eqref{eq:main growth condition} in terms of $E^*_N$, which we recall counts the number of solutions to the inequality 
    \[
    |x_{n_1}-x_{n_2}+x_{n_3}-x_{n_4}| < 1 
    \]
in variables $1\leq n_1,n_2,n_3,n_4\leq N$. For positive integers $k$, define 
    \[
    r(k) = \sum_{\substack{1\leq n_1<n_2\leq N\\k\le x_{n_2}-x_{n_1}< k+1}}1, 
    \]
and denote by $\cR$ the support of $r$; namely, 
    \[
    \cR = \{k\in\NN:r(k)\neq0\}. 
    \]
We note that 
    \begin{equation}
    \label{eq:R size}
    |\cR|\leq N^2, 
    \end{equation}
since in the extremal case where each difference $x_{n_2}-x_{n_1}$ falls into a distinct interval of the form $(k,k+1]$, there can be at most $N^2$ values of $k$ with $r(k)\neq 0$. We have 
    \[
    \sum_{\substack{1\leq n_1<n_2\leq N\\x_{n_2}-x_{n_1}\geq1}}\frac{1}{(x_{n_2}-x_{n_1})^{1/2}} = \sum_{k=1}^\infty\sum_{\substack{1\leq n_1<n_2\leq N\\k\leq x_{n_2}-x_{n_1}<k+1}}\frac{1}{(x_{n_2}-x_{n_1})^{1/2}} \ll \sum_{k\in\cR}\frac{r(k)}{k^{1/2}}, 
    \]
and by the Cauchy--Schwarz inequality, 
    \[
    \sum_{k\in\cR}\frac{r(k)}{k^{1/2}} \leq \left(\sum_{k\in\cR}\frac{1}{k}\right)^{1/2}\left(\sum_{k\in\cR}r(k)^2\right)^{1/2}. 
    \]
Using~\eqref{eq:R size} gives
    \[
    \sum_{k\in\cR}\frac{1}{k} \leq \sum_{\ell\leq N^2}\frac{1}{\ell} \ll N^{o(1)}, 
    \]
and for summation over $r(k)$, we have 
    \begin{align*}
    \sum_{k\in \cR}r(k)^2 &= \sum_{k\in\cR}\sum_{\substack{1\leq n_1<n_2\leq N\\k\le x_{n_2}-x_{n_1}< k+1}}\sum_{\substack{1\leq n_3<n_4\leq N\\k\le x_{n_4}-x_{n_3}< k+1}}1 \\
    &\ll \sum_{\substack{1\leq n_1<n_2\leq N\\1\leq n_3<n_4\leq N\\|x_{n_2}-x_{n_1}-x_{n_4}+x_{n_3}|<1}}1 \\
    &\ll E^*_N. 
    \end{align*}
It follows that 
    \[
    \sum_{k\in\cR}\frac{r(k)}{k^{1/2}} \ll N^{o(1)}(E^*_N)^{1/2}
    \]
which implies that 
    \[
    \sum_{\substack{1\leq n_1<n_2\leq N \\ x_{n_2}-x_{n_1}\ge 1}}\frac{1}{(x_{n_2}-x_{n_1})^{1/2}} \ll N^{o(1)}(E^*_N)^{1/2}. 
    \]
This completes the proof upon noting that the conditions~\eqref{eq:main energy condition} and~\eqref{eq:main growth condition} are satisfied with $\kappa=1/2$ under the assumption~\eqref{eq:energy only condition}. 

\subsection{Concentration in short intervals: proof of Corollary~\ref{cor:short interval}}

Applying Theorem~\ref{thm:main} with $\delta\to\delta/2$ and $\kappa=\delta/2$, it is sufficient to show that 
    \[
    \sum_{1\leq n_1<n_2\leq N}\frac{1}{(x_{n_2}-x_{n_1})^{1/2}} \ll N^{1+o(1)}. 
    \]
Recall the assumption~\eqref{eq:growth condition short interval}
    \[
   \sum_{\substack{n\in\NN\\x_n\in[X,X+H]}}1 \ll H^{1/2}. 
    \]
We have 
    \begin{align*}
    \sum_{\substack{1\leq n_1<n_2\leq N\\x_{n_2}-x_{n_1}\geq1}}\frac{1}{(x_{n_2}-x_{n_1})^{1/2}} &\ll \sum_{1\le n_1\le N}\sum_{j=0}^\infty\frac{1}{2^{j/2}}\sum_{\substack{1\leq n_2\leq N\\x_{n_2}\in[x_{n_1}+2^j,x_{n_1}+2^{j+1}]}}1 \\
    &\ll \sum_{1\leq n_1\leq N}\sum_{j=0}^\infty\frac{\min\{N,2^{j/2}\}}{2^{j/2}} \\
    &\ll N^{1+o(1)}, 
    \end{align*}
which completes the proof. 

\section{Energy estimates}
\label{sec:energy}

In the light of Corollary~\ref{cor:energy only}, the problem of determining whether a sequence is metric Poissonian can be approached via establishing a satisfactory energy estimate $E^*_N$, defined by 
    \[
    E^*_N = |\{1\leq n_1,n_2,n_3,n_4\leq N:|x_{n_1}-x_{n_2}+x_{n_3}-x_{n_4}|<1\}|. 
    \]
The main result for quantitatively convex sequences is as follows. 

\begin{proposition}
\label{prop:convex}
Let $(x_n)_{n\in\NN}$ be a sequence of real numbers satisfying 
    \begin{equation}
    \label{eq:convex condition again}
    x_{n+1}-x_n \geq x_{n}-x_{n-1}+cn^{-1/10^{4}}. 
    \end{equation}
Then 
    \[
    E^{*}_N\ll N^{247/100+o(1)}. 
    \]
\end{proposition}

One could establish the result under a weaker assumption by sharpening the exponent $10^{-4}$ in~\eqref{eq:convex condition again}. We have made no attempt to do so, and for this reason, we have also not tried to optimise the exponent in the bound for $E^*_N$. 

A polynomial $p\in\RR[X]$ of degree at least $2$ with positive leading coefficient satisfies 
    \[
    p(n+1)-p(n) \geq p(n)-p(n-1)
    \]
for all sufficiently large $n$, and therefore also satisfies the condition~\eqref{eq:convex condition again}. Notably, this covers all cases relevant to the metric Poissonian problem for polynomial sequence, since we may assume sequences are eventually increasing. 

The following is a separate result for polynomial sequences that does not require the positivity assumption on the leading coefficient, which also benefits from an alternative approach that avoids the use of advance techniques from additive combinatorics. 

\begin{proposition}
\label{prop:polynomial}
For any polynomial $p\in\RR[X]$ of degree $k\geq2$, 
    \[
    E^*_N \ll N^{3-2/(k+2)+o(1)}.
    \]
\end{proposition}

Now Corollary~\ref{cor:short interval} applies to conclude that non-linear polynomial sequences are metric Poissonian, since the condition~\eqref{eq:growth condition short interval} is satisfied; namely, for any polynomial $p$ of degree $k\geq2$ and $X,H\gg 1$, 
    \[
    |\{n\in\NN:p(n)\in [X,X+H]\}|\ll H^{1/k}. 
    \]
Note also that quadratic polynomials satisfy the stronger estimate $E^*_N\ll N^{2+o(1)}$ (see~\cite[Lemma~5.2]{BKW10}), which is covered by Theorem~\ref{thm:AEM} of Aistleitner, El-Baz and Munsch. 

The preliminaries and proofs of Propositions~\ref{prop:convex} and~\ref{prop:polynomial} are given in Sections~\ref{sec:convex} and~\ref{sec:polynomial}, respectively, thereby proving Corollaries~\ref{cor:convex sequences} and~\ref{cor:polynomial sequences}. 

\subsection{Convex sequences: proof of Corollary~\ref{cor:convex sequences}}
\label{sec:convex}

Let $\cA$ be a finite set of real numbers and $E(\cA)$ denote the additive energy of $\cA$, defined by 
    \[
    E(\cA) = |\{a_1,a_2,a_3,a_4\in\cA:a_1+a_2=a_3+a_4\}|. 
    \]
We employ a discretisation argument in order to appeal to the following result of Bloom~\cite[Theorem~2]{B25+}. 

\begin{lemma}
\label{lem:bloom}
Let $\cA=\{a_1,\dots,a_N\}$ be a set of real numbers satisfying 
    \[
    a_{n+1}-a_{n} > a_n-a_{n-1} \quad \text{for} \quad 1<n<N. 
    \]
Then 
    \[
    E(\cA) \ll |\cA|^{123/50+o(1)}. 
    \]
\end{lemma}

Our reduction to Lemma~\ref{lem:bloom} requires bounds on how the additive energy behaves with respect to unions of sets. Results of this sort are well-known and we provide a proof of the following lemma based on Minkowski's inequality. 

\begin{lemma}
\label{lem:gowers}
Let $\cA_1,\dots,\cA_K$ be disjoint sets of integers. Then 
    \[
    E\left(\bigcup_{j=1}^{K}\cA_j\right)^{1/4} \ll \sum_{j=1}^{K}E(\cA_j)^{1/4}. 
    \]
\end{lemma}
\begin{proof}
Let $\1_{\cA_j}$ denote the indicator function of the set $\cA_j$. Then 
    \[
    E\left(\bigcup_{j=1}^{K}\cA_j \right)^{1/4} = \sum_{j=1}^{K}\left(\int_{0}^{1}\left|\sum_{n\in\ZZ}\1_{\cA_j}(n)e(nx)\right|^{4}dx\right)^{1/4}. 
    \]
By Minkowski's inequality, 
    \begin{align*}
    E\left(\bigcup_{j=1}^{K}\cA_j \right)^{1/4} &\leq \left(\int_{0}^{1}\left|\sum_{j=1}^{K}\sum_{n\in\ZZ}\1_{\cA_j}(n)e(nx)\right|^{4}dx\right)^{1/4} \\
    &= \sum_{j=1}^KE(\cA_j)^{1/4}. \qedhere
    \end{align*}
\end{proof}

\begin{proof}[Proof of Proposition~\ref{prop:convex}]
Let 
    \begin{equation}
    \label{eq:Kdef1}
    K = \left\lfloor N^{1/100}\right\rfloor 
    \end{equation}
and for each $n\leq N$, we define 
    \[
    X_n = \lfloor Kx_n\rfloor. 
    \]
The condition~\eqref{eq:convex condition again} implies that for any $n\in\NN$, 
    \begin{equation}
    \label{eq:convex condition weak}
    x_{n+1}-x_n \gg n^{1-1/10^4}, 
    \end{equation}
and so $X_1,\dots,X_N$ are distinct integers. Moreover, we have  
    \begin{equation}
    \label{eq:AAAA}
    \left|\frac{X_n}{K}-x_n\right| \ll \frac{1}{K}. 
    \end{equation}
If $n_1,n_2,n_3,n_4$ satisfy 
    \[
    |x_{n_1}-x_{n_2}+x_{n_3}-x_{n_4}| \leq 1, 
    \]
then 
    \[
    |X_{n_1}-X_{n_2}+X_{n_3}-X_{n_4}| \ll K|x_{n_1}-x_{n_2}+x_{n_3}-x_{n_4}|+1 \ll K.
    \]
It follows that for some absolute constant $C$, 
    \begin{equation}
    \label{eq:energy convex S(j) bound}
    E^*_N \leq \sum_{\substack{1\leq n_1,n_2,n_3,n_4\leq N\\|X_{n_1}-X_{n_2}+X_{n_3}-X_{n_4}|\leq CK}}1 = \sum_{j=-CK}^{CK}S(j)
    \end{equation}
where 
    \[
    S(j) := \sum_{\substack{1\leq n_1,n_2,n_3,n_4\leq N\\X_{n_1}-X_{n_2}+X_{n_3}-X_{n_4}=j}}1. 
    \]
By orthogonality of additive characters and the triangle inequality, 
    \begin{align*}
    S(j) &= \sum_{1\leq n_1,n_2,n_3,n_4\leq N}\int_0^1e((X_{n_1}-X_{n_2}+X_{n_3}-X_{n_4}-j)t)dt \\
    &= \int_0^1\left(\sum_{1\le n \le N}e(X_nt)\right)^2\left(\sum_{1\le n \le N}e(-X_nt)\right)^2e(-jt)dt \\
   & \le \int_0^1\left|\sum_{1\le n \le N}e(X_nt)\right|^4dx \\
   &= S(0). 
    \end{align*}
Hence, by~\eqref{eq:Kdef1} and~\eqref{eq:energy convex S(j) bound}, 
    \begin{equation}
    \label{eq:E*}
    E^*_N \ll KS(0) \leq N^{1/100}E(\cX) 
    \end{equation}
where 
    \[
    \cX := \{X_j:1\leq j\leq N\}. 
    \]
We partition $\cX$ into subsets 
    \[
    \cX_k = \{X_j\in\cX:j\equiv k\Mod{K}\}, 
    \]
and apply Lemma~\ref{lem:gowers} to get 
    \begin{equation}
    \label{eq:XA}
    E(\cX)^{1/4} \leq \sum_{k=1}^KE(\cX_k)^{1/4}. 
    \end{equation}

We next show the sets $\cX_k$ are convex. Let 
    \[
    X_{(j-1)K+k}, \quad X_{jK+k} \quad \text{and} \quad X_{(j+1)K+k}
    \]
be three consecutive elements of $\cX_k$. By~\eqref{eq:AAAA}, there exists a constant $C$ such that 
    \begin{align*}
    X_{(j+1)K+k}-X_{jK+k} &\geq K(x_{(j+1)K+k}-x_{jK+k})-C \\
    &=K\sum_{\ell=1}^{K}(x_{jK+k+\ell}-x_{jK+k+\ell-1})-C.
    \end{align*}
For $N$ sufficiently large, by~\eqref{eq:convex condition weak}, we have 
    \[
    x_{jK+k+\ell}-x_{jK+k+\ell-1} \geq x_{(j-1)K+k+\ell}-x_{(j-1)K+k+\ell-1}+N^{1/400}.
    \]
This implies 
    \begin{align*}
    X_{(j+1)K+k}-X_{jK+k} &\geq K\sum_{\ell=1}^{K}x_{(j-1)K+k+\ell}-x_{(j-1)K+k+\ell-1}+N^{1/800} \\
    &=K(x_{jK+k}-x_{(j-1)K+k})+N^{1/800}.
    \end{align*}
Another application of~\eqref{eq:AAAA} results in 
    \begin{align*}
    X_{K(j+1)+k}-X_{Kj+k} &\geq X_{jK+k}-X_{(j-1)K+k}+N^{1/900} \\ &> X_{Kj+k}-X_{K(j-1)+k}, 
    \end{align*}
again assuming $N$ is sufficiently large. This implies that the set $\cX_k$ is convex and hence by Lemma~\ref{lem:bloom}, 
    \[
    E(\cX_k) \ll |\cX_k|^{123/50+o(1)}. 
    \]
Finally, by~\eqref{eq:XA}, 
    \[
    E(\cX)^{1/4} \ll N^{o(1)}\sum_{k=1}^K|\cX_k|^{123/200} \ll N^{123/200+o(1)}K^{1-123/200}, 
    \]
which after combining with~\eqref{eq:Kdef1} and~\eqref{eq:E*} gives
    \[
    E^{*}_N\ll N^{1/100}N^{123/50+o(1)}\le N^{247/100+o(1)}. \qedhere
    \]
\end{proof}

\subsection{Polynomial sequences: proof of Corollary~\ref{cor:polynomial sequences}}
\label{sec:polynomial}

We may assume $k\geq3$. Our proof is based on a generalisation of Hua's inequality and is similar to some arguments occurring in~\cite[Section~14]{W19}. 

We require the following result of Bourgain, Demeter and Guth~\cite[Theorem~4.1]{BDG16} and of Wooley~\cite[Theorem~1.1]{W19} on the Vinogradov mean value. 

\begin{lemma}
\label{lem:vmvt}
Let $k\ge 3$ be an integer, $s\leq k(k+1)/2$ be a positive real number, and $a_1,\dots,a_N$ be complex numbers. Then 
    \[
    \int_{[0,1)^{k}} \left|\sum_{1\le n \le N}a_n e(y_1 n+\ldots+y_k n^{k}) \right|^{2s}dy_1\ldots dy_{k} \ll H^{o(1)}\left(\sum_{1\le n \le N }|a_n|^2 \right)^{s}.
    \]
\end{lemma}

Our next result follows from combining Lemma~\ref{lem:vmvt} with some ideas from the the traditional proof of Hua's inequality.

\begin{lemma}
\label{lem:hua}
For any polynomial $p\in\RR[X]$ of degree $k\geq3$, 
    \[
    \int_0^1\left|\sum_{1\leq n\leq N}e(p(n)x)\right|^{k(k+1)}dx \ll N^{k^2+o(1)}. 
    \]
\end{lemma}

\begin{proof}
Let 
    \[
    I = \int_0^1\left|\sum_{1\leq n\leq N}e(p(n)x)\right|^{k(k+1)}dx.
    \]
Applying the identity 
    \[
    1 = 4N^{j}\alpha_j(n)\int_{-1/(4N^{j})}^{1/(4N^{j})}e(tn^{j})dt 
    \]
where 
    \[
    \alpha_j(n) = \frac{\pi n^{j}/(4N^j)}{\sin(2\pi n^{j}/(4N^{j}))}, 
    \]
we obtain 
    \begin{align*}
    I &= (4N)^{k^2(k-1)(k+1)/2} \\
    &\quad \times\int_0^1\left|\int_{B_N}\sum_{1\leq n\leq N}\alpha(n)e\left(p(n)x+\sum_{j=1}^{k-1}x_jn^j\right)dx_1\dots dx_{k-1}\right|^{k(k+1)}dx 
    \end{align*}
where 
    \[
    \alpha(n) := \prod_{j=1}^{k-1}\alpha_j(n) 
    \]
and $B_N$ denotes the box 
    \[
    B_N = \prod_{j=1}^{k-1}\left[-\frac{1}{4N^j},\frac{1}{4N^{j}}\right].
    \]
By H\"{o}lder's inequality, 
    \begin{align*}
    I &\ll N^{k(k-1)/2} \\
    &\quad \times\int_{[0,1]\times B_N}\left|\sum_{1\leq n\leq N}\alpha(n)e(p(n)x+x_{k-1}n^{k-1}+\dots+x_1n)\right|^{k(k+1)}d\tilde{x}
    \end{align*}
where $d\tilde{x}:=dx_1\dots dx_{k-1}dx$. Since the polynomial has degree $k$, we may write 
    \[
    p(n) = c_kn^k+\cdots+c_1n+c_0
    \]
with $c_k\neq 0$, and so 
    \begin{align*}
    I &\ll N^{k(k-1)/2}\int_{[0,1]\times B_N} \left|\sum_{1\leq n\leq N}a(n)e\left(c_kxn^k+\sum_{j=1}^{k-1}(c_jx+x_j)n^j\right)\right|^{k(k+1)}d\tilde{x}.
    \end{align*}
Applying the change of variables 
    \[
    y_j = c_jx+x_j \quad \text{and} \quad y = c_kx, 
    \] 
and noting that the corresponding Jacobian matrix is given by 
    \[
    J := \begin{pmatrix}
    1 & 0 & \cdots & 0 & c_1 \\
    0 & 1 & \cdots & 0 & c_2 \\
    \vdots & \vdots & \ddots & \vdots & \vdots \\
    0 & 0 & \cdots & 1 & c_{k-1} \\
    0 & 0 & \cdots & 0 & c_k
    \end{pmatrix}
    \]
and so $\det J=c_k$, we obtain 
    \begin{align*}
    I &\ll N^{k(k-1)/2}\int_{[-C,C]^{k}} \left|\sum_{1\leq n\leq N}a(n)e\left(yn^k+\sum_{j=1}^{k-1}y_jn^j\right)\right|^{k(k+1)}d\tilde{y}
    \end{align*}
where $d\tilde{y}:=dy_{1}\dots dy_{k-1}dy$ and $C$ is a sufficiently large constant depending only on the polynomial $p$. By periodicity, 
    \begin{align*}
    I &\ll N^{k(k-1)/2}\int_{[-1,1]^{k}} \left|\sum_{1\leq n\leq N}a(n)e\left(yn^k+\sum_{j=1}^{k-1}y_jn^j\right)\right|^{k(k+1)}d\tilde{y}.
    \end{align*}
Combining Lemma~\ref{lem:vmvt} with the fact that $\alpha(n)\ll1$, we obtain
    \[
    I \ll N^{k(k-1)/2}\cdot N^{k(k+1)/2+o(1)} = N^{k^2+o(1)}. \qedhere
    \]
\end{proof}

\begin{proof}[Proof of Proposition~\ref{prop:polynomial}]
Recall that 
    \[
    E_N^* = |\{1\leq n_1,n_2,n_3,n_4\leq N:|p(n_1)-p(n_2)+p(n_3)-p(n_4)|<1\}|. 
    \]
By Lemma~\ref{lem:mean value theorem}, 
    \begin{equation}
    \label{eq:Sub!}
    E_N^* \ll \int_{0}^{1}\left|\sum_{1\le n\le N}e(xp(n))\right|^{4}dx. 
    \end{equation}
Let 
    \begin{equation}
    \label{eq:alpha exponent def}
    \alpha(k) = \frac{4}{k(k+1)-2}, 
    \end{equation}
and note that the assumption $k\ge 3$ implies 
    \[
    \alpha(k) \leq \frac{2}{5}. 
    \]
We write 
    \begin{align*}
    &\quad \int_{0}^{1}\left|\sum_{1\le n\le N}e(xp(n))\right|^{4}dx \\
    &= \int_{0}^{1}\left|\sum_{1\le n\le N}e(xp(n))\right|^{2-\alpha(k)}\left|\sum_{1\le n\le N}e(xp(n))\right|^{2+\alpha(k)}dx, 
    \end{align*}
and apply H\"{o}lder's inequality to get 
    \begin{align*}
    E_N^* &\ll \left(\int_{0}^{1}\left|\sum_{1\le n\le N}e(xp(n))\right|^{2}dx\right)^{1-\alpha(k)/2} \\
    &\quad \times\left(\int_{0}^{1}\left|\sum_{1\le n\le N}e(xp(n))\right|^{k(k+1)}dx\right)^{\alpha(k)/2}. 
    \end{align*}
By Lemma~\ref{lem:mean value theorem}, 
    \[
    \int_{0}^{1}\left|\sum_{1\le n\le N}e(xp(n))\right|^{2}dx \ll \sum_{\substack{1\leq n_1,n_2\leq N\\|p(n_1)-p(n_2)|\leq1}}1 \ll N, 
    \]
and by Lemma~\ref{lem:hua}, 
    \[
    \int_{0}^{1}\left|\sum_{1\le n\le N}e(xp(n))\right|^{k(k+1)}dx \ll N^{k^2+o(1)}.
    \]
Substituting the above and using~\eqref{eq:alpha exponent def}, we obtain 
    \[
    E_N^* \ll N^{1-\alpha(k)/2+k^2\alpha(k)/2+o(1)}\ll N^{3-2/(k+2)+o(1)}. \qedhere
    \]
\end{proof}

\bibliographystyle{amsplain}
\bibliography{references}

\end{document}